\documentclass[10pt,reqno]{amsart} 
\usepackage{center-preamble}
\usepackage{etoolbox}
\title[The center of the REA via quantum minors]{The center of the reflection equation algebra via quantum minors}
\author{David Jordan and Noah White}

\begin{document}

\begin{abstract}
We give simple formulas for the elements $c_k$ appearing in a quantum Cayley-Hamilton formula
for the reflection equation algebra (REA) associated to the quantum group $U_q(\gl)$, answering a question of Kolb and Stokman.  The $c_k$'s are certain canonical generators of the center of the REA, and hence of $U_q(\gl)$ itself; they have been described by Reshetikhin using graphical calculus, by Nazarov-Tarasov using quantum Yangians, and by Gurevich, Pyatov and Saponov using quantum Schur functions; however no explicit formulas for these elements were previously known.

As byproducts, we prove a quantum Girard-Newton identity relating the $c_k$'s to the so-called quantum power traces, and we give a new presentation for the quantum group $U_q(\gl)$, as a localization of the REA along certain principal minors.
\end{abstract}
\maketitle
\setcounter{tocdepth}{2}
\tableofcontents

\section{Introduction}
\label{sec:introduction}

Attached to an $N\times N$ matrix $A$ is its \emph{characteristic polynomial},
$$p_A(t) =  t^N - c_1\cdot t^{N-1} + \cdots + (-1)^N c_N.$$
Here $c_1$ and $c_N$ are the trace and determinant of $A$; more generally each $c_k$ is  the sum of all principal $k\times k$ minors of $A$, and can be expressed as the $k$th elementary symmetric function in the eigenvalues of $A$.  Taken together, the functions $c_1,\ldots, c_N$ generate the algebra $\mathcal{O}(Mat_N)^{GL_N}$ of conjugation-invariant polynomial functions on the variety of $N\times N$ matrices. 
The Cayley-Hamilton theorem states that $p_A(A)= 0$ identically, i.e.
\begin{equation}
  \label{eqn:CH}
  A^N - c_1\cdot A^{N-1} + \cdots + (-1)^N c_N = 0.
\end{equation}
In addition to the coefficients $c_k$ of the characteristic polynomial, there are the power traces $s_k=tr(A^k)$, for $k=1,\ldots, N$; these can be expressed as the power sum symmetric functions of the eigenvalues of $A$.  The elementary and power sum generators are related via the Girard-Newton identities,
\begin{equation} 
  \label{eqn:NI}
  kc_k = \sum_{i=1}^k (-1)^{i-1}s_ic_{k-i},
\end{equation}
which allow us to write either set of generators uniquely in terms of the others.

The present paper is concerned with the center of the \emph{reflection equation algebra} (REA), denoted $\REA$, associated to the quantized universal enveloping algebra $U_q(\gl)$, and with the analog of identities \eqref{eqn:CH} and \eqref{eqn:NI}.   The REA is generated by the entries, $a^i_j$, of an $N\times N$ matrix; the commutation relations amongst the $a^i_j$ are encoded in the so-called reflection equation relations, see Equation \eqref{eqn-RE}. The REA defines a flat (PBW) $q$-deformation of the coordinate algebra $\mathcal{O}(Mat_N)$, along the \emph{Semenov-Tian-Shansky} Poisson bracket \cite{M06}.

The reflection equation first arose in the context of mathematical physics, specifically integrable systems related to factorizable scattering on a half-line with reflecting wall. The RE prescribes the consistency relations for collisions with the reflecting wall, analogously to the QYBE for collisions on the interior (see \cite{K96,S92,C84,KS93}).  Somewhat independently, the REA played a central role in the braided Hopf algebra perspective on quantum groups pioneered by Majid \cite{M91,LM94,M95} and the quantization of conjugacy classes by Donin and Mudrov \cite{DM02, M07}.  More recently, the REA associated to $U_q(\gl)$ arose as the Hochschild homology/horizontal trace of the category $\Rep_q(GL_N)$ of representations of the quantum group \cite{BZBJ15,BZBJ16}, and in the formulation of a $K$-theoretic geometric Satake theorem \cite{CK15}. Finally, we note that the reflection equation appears in the context of co-ideal subalgebras \cite{N96,MRS03}, in particular the reflection equation algebra is an important source of co-ideal sub-algebras admitting a universal $K$-matrix \cite{KS09,K17,JM11} making them into braided module categories for $\Rep_q(GL_N)$.

The conjugation action of $GL_N$ on $Mat_N$ also $q$-deforms to the REA, making it into a $U_q(\gl)$-equivariant algebra.  The Joseph-Letzter-Rosso isomorphism \cite{JL92, R90thesis, R90} gives an algebra embedding of the REA into the locally finite sub-algebra of $U_q(\gl)$.  When $q$ is not a root of unity, the so-called Harish-Chandra isomorphism proved in \cite{JL92} (see also \cite[Proposition 3.6]{KS09}) identifies the center of the REA with its invariant subalgebra; this in turn is isomorphic abstractly to a polynomial ring in $N$ commuting generators.

While these isomorphisms describe the center of $U_q(\gl)$ and $\REA$ abstractly, they do not yield any concrete formulas for generators, nor any quantum analog of the Cayley-Hamilton and Girard-Newton identities.  The combinatorial problem of giving such formulas has a long history.  In Section \ref{sec-related} we review: the diagrammatic techniques of Reshetikhin \cite{R89,B98}, the $q$-Schur function approach of \cite{GPS97,PS96,IOP99}, the quantum Yangian approach of \cite{NT94}, and finally the parallel works of \cite{DL03} and \cite{Z98} for the FRT algebra.

However, what is still missing from all these developments is any closed-form, PBW-ordered expression for a system of generators of the center of the REA, and in particular for those appearing as coefficients in a quantum Cayley-Hamilton identity.  Remarkably, there is to date no such expression even for the quantum determinant in $\REA$ (except for $N=2$, see Example~\ref{exm:ck-n=2and3}), which must appear when giving a generators-and-relations presentation of $\REAG$ and $\mathcal{O}_q(SL_N)$.  This issue has been highlighted in \cite{KS09}, Remark 5.10, and \cite{J14}, Remark 7.9 (see also \cite[Section~4]{DL05} for a related question about explicit computations in the REA). This state of affairs is rectified in this paper: we produce simple, explicit and PBW-ordered expressions for each of the central generators $c_1,\ldots c_N$ appearing in a quantum Cayley-Hamilton identity, and we identify our elements with the prior formulations given in Section \ref{sec-related}.

\subsection{Statement of Results}
\label{sec:statement-results}
We will need the following notation:  for each $N\in\mathbb{N}$, we have the set $\SetN = \{1,\ldots,N\}$, we have the collection ${\SetN \choose k}$ of its $k$-element subsets, and we have the symmetric group $S_N$ of permutations of $\SetN$.  For each $I\in {\SetN \choose k}$, we have the subgroup $\Sym(I)\subset S_N$, consisting of all elements which fix the complement $I^\comp$ of $I$ in $\SetN$, and we have the \emph{weight}, $\wt(I) = \sum_{i\in I} i.$  Finally, for each permutation $\sigma\in S_N$, we have its \emph{length} and its \emph{exceedance}:
\begin{align*}\lng(\sigma) &= \textrm{ the number of pairs $i<j\in\SetN\times \SetN$ such that $\sigma(i)>\sigma(j)$},\\
e(\sigma) &= \textrm{ the number of $i\in \SetN$ such that $\sigma(i)>i$}.\end{align*}

For a set \( I \in {[N]\choose k} \) we will index its elements in ascending order \( i_1 < i_2 < \ldots < i_k \). Finally, recall that the generators of $\REA$ organize naturally into the $N\times N$ matrix,
$$A = \sum_{i,j}a^i_jE^j_i.$$
In the following discussion we regard $A$ as an element in $Mat_N(\REA)$. With this notation in hand, we can state our main result:
\begin{Theorem}[see Corollaries~\ref{cor:cs-are-DL-invars} and~\ref{cor:qCH}]\label{main-thm}
The center of $\mathcal{O}_q(Mat_N)$ is freely generated by the elements:
$$c_k = \!\!\sum_{I \in {\SetN \choose k}}q^{-2\wt(I)}\!\!\sum_{\sigma\in Sym(I)} (-q)^{\ell(\sigma)}\cdot q^{e(\sigma)}\cdot a^{i_1}_{\sigma(i_1)}\cdots a^{i_k}_{\sigma(i_k)},$$
for $k=1,\ldots, N$.  Moreover, the elements $c_k$ satisfy the following \emph{quantum Cayley Hamilton} identity:
\begin{equation}
  \sum_{k=0}^N(-q^2)^{N-k}c_{N-k}\cdot A^k =0 \label{eqn:QCH}.
\end{equation}
\end{Theorem}

\begin{Example}
  \label{exm:ck-n=2and3}
In the case $N=2$, our formulas (up to a choice of normalization) appeared already in \cite{M91} and \cite{KS92}:
\begin{equation*}
  c_1 = q^{-2}a^1_1 + q^{-4} a^2_2, \qquad c_2  = q^{-6}\left(a^1_1a^2_2 - q^2a^1_2a^2_1\right).
\end{equation*}
Already for $N=3$, our formulas for $c_2$ and $c_3$ appear to be new. They read:
\begin{align*}
  c_1 &= q^{-2}a^1_1+q^{-4}a^2_2 + q^{-6}a^3_3\\
  c_2 &= q^{-6}(a^1_1a^2_2-q^2a^1_2a^2_1) + q^{-8}(a^1_1a^3_3 - q^4a^1_3a^3_1) + q^{-10} ( a^2_2a^3_3 - q^2a^2_3a^3_2),\\
  c_3 &=  q^{-12}\left(a^1_1a^2_2a^3_3 - q^2 a^1_1a^2_3a^3_2 -q^2 a^1_2a^2_1a^3_3 -q^4 a^1_3a^2_2a^3_1+ q^4 a^1_2a^2_3a^3_1 +q^3a^1_3a^2_1a^3_2\right).
\end{align*}
These should be contrasted with superficially similar formulas for the conjugation co-invariants for the FRT algebra given in Example \ref{DLinv-n=2and3}.
\end{Example}

\begin{Remark} Note that we have an obvious isomorphism $\Sym(I) \cong S_{|I|}$, however the notion of length appearing in the formula above is that inherited from $S_N$, not from $S_{|I|}$. The first instance of this is the coefficient $q^4$ of $a^1_3a^3_1$  above.
\end{Remark}

\subsection{Relations to other works}\label{sec-related}

Our formulas for $c_k$ are motivated by a number of well-known descriptions for the center of the REA and/or $U_q(\gl)$.  Let us now outline these antecedents.

Our starting point is in fact Reshetikhin's description \cite{R89} (see also \cite{B98}), which uses the graphical calculus for the REA, as championed independently by Lyubashenko-Majid \cite{LM94}.  Reshetikhin's canonical generators $c_k$ are defined as the ``quantum traces'' of the irreducible representations $V(\omega_i)$, for $i=1,\ldots N$.  Here, by quantum trace, we mean the composite,
$$c_k:\mathbf{1} \xrightarrow{\coev} V(\omega_i)\otimes V(\omega_i)^* \xrightarrow{\sigma} V(\omega_i)^* \otimes V(\omega_i)\hookrightarrow \REA.$$
See Figure~\ref{fig-ckid} for a depiction of $c_k$ in the graphical calculus, and Section~\ref{sec:graph-calc} for a discussion of the graphical caluclus for $\REA$ more generally.

The explicit formulas we produce in Theorem~\ref{main-thm} are more precisely formulas for Reshetikhin's central elements $c_k$.  The essential challenge here is that while Reshetikhin's generators have a simple and elegant formulation in the 3-dimensional graphical calculus, translating these into ``1-dimensional'' (i.e, algebraic) formulas is highly non-trivial.

The first appearance of a quantum Cayley-Hamilton identity in the context of quantum groups appears in the work of Nazarov and Tarasov \cite{NT94}.  Their construction involves realizing $U_q(\gl)$ as a quotient of the quantum Yangian,
and establishing the required identities in the quantum Yangian.

The quantum Cayley-Hamtilon identity was systematically revisited in work of Gurevich, Pyatov, and Saponov.  In \cite{GPS97} two systems of central generators $s_k$ and $\sigma_k$ are proposed as $q$-deformations of Schur functions, namely the power sum and elementary symmetric functions respectively.  Perhaps the simpler of these are the quantum power traces, $s_k=tr_q(A^k)$, which are $q$-deformations of the classical power traces.  Here $tr_q$ of a matrix $B$ is defined as the weighted sum $\sum_{i}b^i_iq^{-2i}$ of its diagonal entries.  The $\sigma_k$ are introduced so as to satisfy a quantum Newton identity with respect to the $s_k$'s.  Such an identity uniquely determines the $\sigma_k$, because the $s_k$'s freely generate the center.  However, they do not propose explicit formulas for the generators $\sigma_k$, such as those of Theorem \ref{main-thm}.  In the same paper, the authors further establish a quantum Cayley-Hamilton identity, precisely \eqref{eqn:QCH} up to matching conventions for scalar factors.

In Theorem~\ref{thm:qNewton-id}, we prove a similar ``quantum Newton identity'':
\begin{equation}
  [k]_q s_k = \sum_{j=1}^{k} \frac{c_{j-1}}{[j-1]_q!}s_{k-j}=0,\label{eqn:QNI}
\end{equation}
relating Reshetikhin's invariants to the quantum power traces.  Our proof uses the graphical calculus and the finite Hecke algebra.  This formula has several consequences.  First, it allows us to identify Reshetikhin's central generators $c_k$ with Gurevich-Pyatov-Saponov's generators $\sigma_k$.  Secondly, it allows us to transport the quantum Cayley Hamilton identity from \cite{GPS97} to Reshetikhin's generators, and hence to see that our expilcit central elements satisfy a quantum Cayley-Hamilton identity.

Finally, let us note that similar, but simpler, formulas to those we obtain for the $c_k$ were discovered in \cite{DL03} for the co-invariants of the so-called FRT algebra (See Example \ref{ex:q-minors}; these feature also in \cite{AY11}).  The FRT algebra is a quantization of $Mat_N$ along the \emph{Sklyanin} Poisson bracket; it becomes isomorphic to the restricted dual Hopf algebra to $U_q(\gl)$, upon inverting its quantum determinant.  The co-invariants in the FRT algebra are not central however (in fact the center is generated by the quantum determinant only), and the analog of the quantum Cayley-Hamilton theorem for them (see \cite{Z98}) does not have a straightforward interpretation as matrix multiplication as for the REA.

There is a process, with many names -- transmutation, equivariantization, twisting -- which produces the REA from the FRT algebra.  The FRT algebra is simpler in many respects than the REA, so we perform computations later in the paper by first identifying Reshetikhin's generators $c_k$ as the twist of Domokos-Lenagan's co-invariants, and then computing the effect of the twisting procedure on these.

We can summarize these relationships as follows:

\begin{Theorem}[see Theorem \ref{thm:qNewton-id}, Corollary \ref{cor:cs-are-DL-invars}]
The central generators $c_k$ asserted in Theorem \ref{main-thm} coincide with:  Reshetikin's central elements $c_k$, the twist of Domokos-Lenagan's central elements $\DLinv{k}$, and (up to a scalar) Gurevich-Pyatov-Saponov's central elements $\sigma_k$.
\end{Theorem}

An interesting and related problem is the \( q \)-analogue of the question asked in \cite[Section~3.2]{MR02}, where Molev and Ragoucy construct central elements of certain reflection equation algebras appearing as subalgebras of the Yangian. 

\subsection{Quantization of conjugacy classes}
The works \cite{DM02, M07} of Donin and Mudrov produce explicit quantizations of semi-simple conjugacy classes in $GL_N$.  These are obtained by quotenting the REA by quantum Cayley-Hamilton relations, specialized at a given central character according to the eigenvalues of the classical orbit.  We hope that the explicit reformulation of the QCH identity in Theorem \ref{main-thm} can further clarify these quantizations and their properties.

As an illustrative example, let us consider the quantization of the unipotent cone in $GL_N$.  These are all matrices whose characteristic polynomial coincides with that of the identity matrix.  For this we consider the co-unit $\epsilon:\REAG\to\CC$, which sends $a^i_j$ to 1 if $i=j$ and to $0$ otherwise.   Restricing to the central elements $c_k$, an immediate computation gives that 
$$\epsilon(c_k) = \sum_{I\in {[N] \choose k}} q^{-2\wt(I)},$$
and hence the resulting quantum Cayley-Hamilton identity becomes:
$$(A-q^{-2})(A-q^{-4})\cdots (A-q^{-2N})=0.$$
Interestingly, while the unipotent cone consists of generically non-semisimple matrices, we see that the defining equation for its natural quantization is in fact multiplicity free.

\subsection{A new presentation of $U_q(\gl)$}
\label{sec:new-pres-u_qgl}
As an application, we can give a new presentation of the quantum group $U_q(\gl)$ as a finite degree extension of a localization of $\REAG$, along certain quantum minors which we introduce herein.  Since the relations of $\REAG$ are completely encoded in the reflection equation algebra, this is entirely independent of the Serre presentation.  It is more closely related to the presentation of $U_q(\gl)$ via so-called $L$-operators (see \cite[Section 8.5]{KS97}).  Since the presentation follows easily from our main result, we outline it here in the introduction.

Let us first recall the Joseph-Letzter-Rosso isomorphism \cite{JL92, R90}, which allows us to identify $\REAG$ with the subalgebra of locally finite vectors inside $U_q(\mathfrak{gl}_N)$ (see also \cite[Proposition 10.34]{KS97} and \cite[Proposition 2.8]{KS09}).  In particular, using this isomorphism, we can just as well regard the central elements $c_k$ in Theorem~\ref{main-thm} (together with the inverse of $c_N$) as generators of the center of $U_q(\mathfrak{gl}_N)$ (in fact this how they were first interpeted by Reshetikhin \cite{R89, B98}).  

Taking this point further, we can obtain a new presentation of $U_q(\mathfrak{gl}_N)$.  We recall the generators $K_i$ in the presentation of $U_q(\mathfrak{gl}_N)$ given in \cite[Section 6.1]{KS97}.  We have the following expression relating $c_N$ to the $K_i$'s:
\begin{equation}c_N = (K_1\ldots K_N)^2.\label{eqn-cNasKs}\end{equation}
It follows\footnote{The statement for $U_q(\mathfrak{gl}_N)$ must be extracted from the corresponding statement for $U_q(\mathfrak{sl}_N)$, but this is straightforward.} from \cite{JL92} (see \cite[Theorem 6.33]{KS97} for an exposition) that while the $q$-commuting element
$$K^{-2\rho} = (K_1^{1-N}K_2^{2-N}\ldots K_{N-1})^2\in U_q(\gl)$$
is a locally finite vector for the quantum adjoint action, the element
$$K^{\rho} = K_1^{N-1}K_2^{N-2}\ldots K_{N-1},$$
is not, and that furthermore all of $U_q(\mathfrak{gl}_N)$ is obtained from its locally finite part by adjoining $K^{\rho}$.

In order to find an expression for $K^\rho$ in the generators of the reflection equation algebra, we need to recall the standard inclusions $U_q(\mathfrak{gl}_{N-1})\hookrightarrow U_q(\mathfrak{gl}_N)$, which can be made compatible with the Joseph-Letzter-Rosso isomorphism.  We encode this as follows.  For $k=1,\ldots, N$ we let $A_{\geq k}$ denote $(N-k+1)\times(N-k+1)$ matrix,
$$A_{\geq k} =  \left(\begin{array}{cccc} a^{k}_{k} & \cdots & a^{k}_N\\ \vdots & \ddots & \vdots\\a^N_{k} & \cdots & a^N_N\end{array}\right).$$
In particular, $A_{\geq 1}=A$, and $A_{\geq N} = (a^N_N)$.

Inspecting~(\ref{eq:REA-explicit-rels}), one sees that the entries of the matrix $A_{\geq k}$ form a subalgebra of $\mathcal{O}_q(GL_N)$ isomorphic to $\mathcal{O}_q(GL_{N-k+1})$.  This is in fact compatible with the inclusion of $U_q(\mathfrak{gl}_{N-k+1})$ into $U_q(\mathfrak{gl}_N)$, upon identifying $\mathcal{O}_q$ in each case with the locally finite vectors.

For an $N\times N$ matrix $A$, we denote by $\detq(A)$ the formula $c_N$ above in the entries of $A$.  The next proposition follows easily from the above considerations and Equation \eqref{eqn-cNasKs}.

\begin{Proposition}\label{dets-and-Ks}
For each $k=1,\ldots, N$, we have:
$$\detq(A_{\geq k}) = (K_{k}\cdots K_N)^2.$$
\end{Proposition}

Now we simply observe that, by Proposition \ref{dets-and-Ks}, we have
$$K^{-2\rho} = \detq(A_{\geq 1})^{(1-N)}\cdot \detq(A_{\geq 2}) \cdots \detq(A_{\geq N}),$$
hence it follows:

\begin{Corollary}
We have an algebra isomorphism,
$$U_q(\gl)\cong \REA[\detq(A_{\geq 1})^{-\frac12}, \ldots, \detq(A_{\geq N})^{-\frac12}].$$
\end{Corollary}

This presentation reflects the following well-known fact:  while $U_q(\gl)$ is often regarded informally as a deformation of the universal enveloping algebra $U(\gl)$, it is more precisely a deformation of the Tits $2^r$-fold cover of the big Bruhat cell $Bw_0B$ inside the group $G$.  The big cell $Bw_0B$ is the complement to the union of hypersurfaces defined by the classical limit of the $\detq(A_{>k})$'s.  Denoting the Tits cover by $\widetilde{Bw_0B}$ and $\REAG^\circ = \REA[\detq(A_{>0})^{-1}, \ldots, \detq(A_{>N-1})^{-1}]$, we summarize the situation with the following diagram:
\[
  \begin{tikzcd}
    \widetilde{Bw_0B} \arrow[two heads,r] \arrow[squiggly,d] & Bw_0B \arrow[hook,r] \arrow[squiggly,d] & G \arrow[squiggly,d] \\
    U_q(\gl) & \REAG^\circ \arrow[hook',l] & \REAG^\circ \arrow[hook',l]
  \end{tikzcd}
\]
where the downward wavy arrows designate ``pass to coordinate algebras and quantize''.

\subsection{Outline of paper}
\label{sec:outline-paper}
In Section~\ref{sec:preliminaries}, we recall the necessary preliminaries about the reflection equation algebra:  its presentation, how to obtain it from the FRT algebra by the twist procedure, and finally we recall the graphical calculus for the REA, in Section~\ref{sec:graph-calc}.  We recall some basic facts about the finite Hecke algebra in~\ref{sec:finite-hecke-algebra}, and use it to introduce Reshetikhin's central generators in~\ref{sec:resh-centr-elem}.  In Section~\ref{sec:quant-girard-newton}, we prove the quantum Girard-Newton identities for Reshetikhin's generators, in the process we relate them to the generators $\sigma_k$ of \cite{PS96}, and hence are able import the quantum Cayley-Hamilton identity from \cite{GPS97} to our setting. In Section~\ref{sec:minors-combinatorics}, we introduce some combinatorics:  we recall in~\ref{sec:permutations} the notions of length and exceedance of a permutation, and then generalize these notions from permutations to general bijections between subsets of a fixed set.  We define the quantum minors which appear in the formulas for $c_k$, and prove a simple row expansion formula for them.  In Section~\ref{sec:expansion-cliques} we introduce the notion of an expansion clique:  this abstracts the set of indices appearing inductively in a row expansion formula for quantum determinants, for use in the proof of our main result.  In Section~\ref{sec:main-theorem}, we prove our main result, i.e. formulas for the central generators $c_k$: this proceeds by an elaborate but elementary induction on the collection of expansion cliques, using the row expansion formula.

\subsection{MAGMA code}
In the process of writing this paper, the authors developed MAGMA \cite{MAGMA} code for working with the reflection equation algebra.  The code is publicly available at \href{http://www.github.com/noaham/REA}{http://www.github.com/noaham/REA}.  The authors would be happy to help an interested reader to use the code.

\subsection{Acknowledgements}
\label{sec:acknowledgements}
We would like to thank Sammy Black, Pavel Etingof and Masahiro Namiki for helpful discussions at an early stage of the project, and Tom Lenagan for helpful comments on a preliminary version of the paper.  Most of the key identities used in the paper were found using extensive computer experimentation in MAGMA.  We thank the Simons Foundation for providing site licenses to US institutions, in particular UCLA, so that we could carry out this experimentation.  The work of the first author is supported by European Research Council (ERC) under the European Union's Horizon 2020 research and innovation programme (grant agreement no. 637618). A portion this work was completed while the second author was a guest at the Max Plank Institute for Mathematics in Bonn.

\section{Preliminaries}
\label{sec:preliminaries}

Here we collect the necessary background on the REA and FRT algebras and the diagrammatic calculus for the REA.  We recall Reshetikhin's central elements \cite{R89}, following the exposition of \cite{B98}.

Let \( R \) denote the quantum R-matrix for the defining representation \( V \) of \( \gl \) (we follow the conventions in \cite{KS97}):

\begin{equation*}\label{eqn-Rmat}
  R = q\sum_{i}E_i^i\otimes E_i^i + \sum_{i\neq j}E_i^i\otimes E_j^j+(q-q^{-1})\sum_{i>j}E_i^j\otimes E_j^i 
\end{equation*}
(an operator on \( V \otimes V \)), and let us denote by $R_{21} = \tau \cdot R \cdot \tau$, where $\tau$ denotes the vector flip $\tau(v\otimes w) = w\otimes v$.  Here $E^i_j$ denotes the elementary matrix, $E^i_j\cdot v_k = v_j,$ if $i=k$, $0$ otherwise. 

For an element \( M \in \End(V\ox V) \), we denote by \( M^{t_2} \in \End(V\ox V) \) the operator obtained by transposing on the second factor. Let \( \tilde{R} = ((R^{t_2})^{-1})^{t_2} \), then, as is easily verified,
\begin{equation*}\label{eqn-Rmat}
  \tilde{R} = q^{-1}\sum_{i}E_i^i\otimes E_i^i + \sum_{i\neq j}E_i^i\otimes E_j^j-(q-q^{-1})\sum_{i>j}q^{2(i-j)}E_i^j\otimes E_j^i.
\end{equation*}

\begin{Definition}[See \cite{KS97}, Example 10.18] The \emph{reflection equation} algebra \( \REA \) has generators \( a^i_j \) for \( i,j \in [N] \).  The relations amongst the generators are the $N^2\times N^2$ entries of the matrix equation,
\begin{equation}\label{eqn-RE}
  R_{21}A_1R_{12}A_2 = A_2R_{21}A_1R_{12},
\end{equation}
where we collect into a matrix \( A = (a^i_j) = \sum a^i_j\otimes E^j_i\), and we denote $A_1=A\otimes\id$, and $A_2=\id\otimes A$.
\end{Definition}
We will need to make use of the explicit form of these relations, which we can write as follows.  Let $\theta$ denote the Heaviside theta function, i.e. \( \theta(x) = 0 \)  if \( x \le 0 \) and \( \theta(x) = 1 \) if \( x > 0 \)  and let \( \delta_{ij} \) denote the Kronecker delta function, so $\delta_{ij}=1$ if $i=j$, zero else.  Applying Equation \eqref{eqn-Rmat}, we have that Equation \eqref{eqn-RE} is equivalent to the list of relations (c.f. \cite[Section~3]{DL05}),
\begin{equation}\label{eq:REA-explicit-rels}
  \begin{aligned}
    &\text{when } j<i:
    &a^i_m a^j_n &= q^{\delta_{i,n} + \delta_{n,m} - \delta_{m,j} } a^j_n a^i_m \\
    &&&\phantom{=} + \theta(n-m) q^{\delta_{i,m} - \delta_{j,m}}(q-q^{-1}) a^j_m a^i_n \\
    &&&\phantom{=} + \delta_{i,n}(q-q^{-1})q^{\delta_{n,m} - \delta_{j,m}} \sum_{p > i} a^j_p a^p_m \\
    &&&\phantom{=} + \delta_{i,m}\theta(n-m)(q-q^{-1})^2 \sum_{p > i} a^j_p a^p_n \\
    &&&\phantom{=} - \delta_{j,m} q^{-1}(q-q^{-1}) \sum_{p > j} a^i_p a^p_n \\
    &\text{when } j=i \text{ and } n<m : \qquad
    &a^i_m a^i_n &= q^{\delta_{i,n} -\delta_{i,m} - 1 } a^i_n a^i_m \\
    &&&\phantom{=} + \delta_{i,n} q^{-1}(q-q^{-1}) \sum_{p > i} a^i_p a^p_m \\
    &&&\phantom{=} - \delta_{i,m} q^{-1}(q-q^{-1}) \sum_{p > i} a^i_p a^p_n.
  \end{aligned}
\end{equation}

\begin{Definition}[See \cite{KS97}, Section 9.1] The FRT algebra, \( \FRT \) has generators \( x^i_j \) for \( i,j \in [N] \). The relations amongst the generators are the $N^2\times N^2$ entries of the matrix equation,
\begin{equation}\label{eqn-FRT}RX_1X_2=X_2X_1R\end{equation}
where we collect into a matrix \( X = (a^i_j) = \sum x^i_j\otimes E^j_i\), and we denote $X_1=X\otimes\id$, and $X_2=\id\otimes X$.
\end{Definition}
Applying Equation \eqref{eqn-Rmat}, we have that Equation \eqref{eqn-FRT} is equivalent to the list of relations,
\begin{equation} \label{eq:FRTrel}
\begin{aligned} 
  x^i_k x^i_l &= q x^i_l x^i_k  \\
  x^i_k x^j_k &= q x^j_k x^i_k  \\
  x^i_l x^j_k &= x^j_k x^i_l  \\
  x^i_k x^j_l - x^j_l x^i_k &= (q-q^{-1}) x^i_l x^j_k 
\end{aligned}
\end{equation}
for any \( 1 \le i < j \le n \) and \( 1 \le k < l \le n \).

\begin{Remark} A key distinction between \( \FRT \) and $\REA$ is that the relations of $\FRT$ are ``local'' in the sense that any \( 2\times 2 \) submatrix of \( X = (x^i_j) \) generates a subalgebra isomorphic to $\widetilde{\mathcal{O}_q}(Mat_2)$.  This is essentially because the $R$ matrix appears only linearly in the relations \eqref{eqn-FRT}, but quadratically in \eqref{eqn-RE}, and it considerably simplifies computations in the FRT algebra versus the REA.
\end{Remark} 

\subsection{The PBW bases for the RE and FRT algebras}
\label{sec:pbw-bases}

Suppose we have an algebra \( A \) given in the form \( A = T(V)/R \), for some finite dimensional vector space \( V \) and a two sided ideal \( R \). A totally ordered basis \( \set{x_1 < x_2 < \ldots < x_r} \) of \( V \) is said to be a \emph{PBW generating set} if the set of ordered monomials
\begin{equation*}
  \setc{ x_{i_1}x_{i_2}\cdots x_{i_k}}{ 1 \le i_1 \le i_2 \le \ldots \le i_k \le r, k \in \mathbb{N}}
\end{equation*}
form a basis of \( A \) which will be called a \emph{PBW basis}.

The generating set \( \set{a^i_j} \) for \( \REA \) is given the lexiographic total order where \( a^i_j < a^k_l \) is \( i < k \) or if \( i=k \) and \( j<l \). The generating set \( \set{x^i_j} \) for \( \FRT \) is given the same total order. For example \( a^2_4 < a^3_2 \) and \( a^2_2 < a^2_3 \) and the monomial \( a^1_3a^3_2a^3_3a^4_1 \) is ordered.

\begin{Proposition}
  \label{prp:PBW-algebras}
  The generating sets \( \set{a^i_j} \) and \( \set{x^i_j} \) are PBW generating sets for the algebras \( \REA \) and \( \FRT \) respectively. 
\end{Proposition}

\begin{proof}
  This is well-known and can be proved directly, see \cite[Proposition~3.1]{DL05} for a careful proof; or alternatively as a consequence of the quantum Peter-Weyl theorem, see for example \cite[Section~11.5.4]{KS97}.
\end{proof}

\subsection{The twist construction}
\label{sec:twist-mult}
We recall here the ``twist'' construction, for relating the algebras $\FRT$ and $\REA$.  We follow the exposition of \cite[Section 10.3]{KS97}, and its categorical reformulation in \cite[Section 3]{JM11}.  For the history of the construction, see \cite{M95}.

To begin, let us recall that the FRT algebra carries commuting $U_q(\gl)$-actions of left and right translation, which are compatible with the co-product in such a way that the left translation action is equivariant for $\Delta^{op}$, rather than $\Delta$.

Hence, the FRT algebra is naturally an algebra in the category
$$\Rep_q(\GL)^{op}\bt\Rep_q(\GL),$$
where $op$ denotes the opposite tensor product, $V\ot^{op} W=W\ot V$.  The tensor product functor,
$$T=(T,\id\bt\sigma\bt\id):\Rep_q(\GL)\bt \Rep_q(\GL)\to \Rep_q(\GL),$$
and the identity functor,
$$\widetilde{Id}=(Id,\sigma): \Rep_q(\GL)^{op}\to\Rep_q(\GL),$$
each are canonically endowed with tensor structures using the braiding in the appropriate place.  Hence we may compose these to obtain a tensor functor,
$$\Rep_q(\GL)^{op}\bt\Rep_q(\GL) \xrightarrow{\widetilde{Id}\bt Id} \Rep_q(\GL)\bt\Rep_q(\GL) \xrightarrow{T} \Rep_q(\GL).$$

Hence, given any algebra in $\Rep_q(\GL)^{op}\bt\Rep_q(\GL)$, its image in $\Rep_q(\GL)$ is again an algebra.  Because the functors $T$ and $Id$ are both compatible with the forgetful functors to vector spaces, this procedure does not change the underlying vector space, but rather it  modifies the multiplication according to the appearances of $\sigma$ in the tensor structures.  It is a straightforward computation with $R$-matrices to see the algebra in $\Rep_q(\GL)$ obtained from $\FRT$ in this way is canonically isomorphic to $\REA$.  As a particular consequence, we have:
\begin{Proposition}[{See \cite[Example 10.18]{KS97} and \cite[Example~1.3.2(c)]{VV10}}] 
  \label{prp:iso}
We have a unique isomorphism of vector spaces,
\[ \Phi\map{\REA}{\FRT}, \] defined on the generators by \[ a^i_j \mapsto x^i_j, \] and extended to arbitrary elements via the cocycle relation,
\begin{equation} \Phi(a b) = m_{FRT}(\mathcal{R}_{13}\mathcal{R}_{23}(a \otimes b)).  \label{eq:mult-relate}\end{equation}
We will denote the inverse of \( \Phi \) by \( \Psi = \Phi^{-1}\map{\FRT}{\REA} \)
\end{Proposition}

\begin{Lemma}[See \cite{KS97}, Example 10.18]
  \label{lem:expanded-rea-mult}
The expression in~(\ref{eq:mult-relate}) is given by
\begin{align*}  \Phi(a^i_j a^k_l) &= R^{in}_{st} x^s_m \tilde{R}^{mk}_{jn} x^t_l \\  &=
    \begin{cases}
      x^i_j x^k_l &\text{if } i < k, j \neq k \\
      q x^i_j x^k_l &\text{if } i=k, j \neq k \\
      q^{-1}x^i_j x^j_l + (q^{-1} - q) \sum_{m > j} q^{-2(m-j)}x^i_m x^m_l &\text{if } i < j, j=k.
    \end{cases}
\end{align*}
and the inverse \( \Psi \) is given explicitly by \( \Psi(x^i_jx^k_l) = (R^{-1})^{ik}_{su}a^s_tR^{tu}_{jv}a^v_l \).
\end{Lemma}

We note in passing that the construction of each of the FRT and RE algebras, as well as the twist construction, can be be understood most conceptually as an instance of the Co-End construction, given originally in this context by Lyubashenko \cite{L95} and Majid \cite{LM94, M95}.  We refer to \cite{JM11} for more details.

\subsection{Graphical calculus and Reshetikhin's central generators}
\label{sec:graph-calc}
In this section we briefly recall the graphical calculus for morphisms in a rigid braided tensor category, referring to \cite{K95} for more details.  We then recall the constructions of Majid which construct the REA as a braided Hopf algebra, with all the structure maps given as simple morphisms in the graphical calculus.  Finally, we use this language to recall a construction Reshetikhin's central generators, adapted to the REA.  We refer to \cite{B98} for a careful exposition of Reshetikhin's elements.

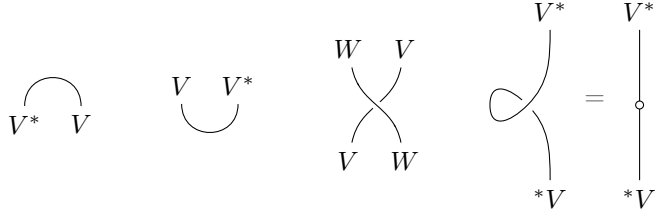
\begin{figure}[h]
  \centering
  \begin{tikzpicture}[xscale=0.75,baseline={(current bounding box.center)}]
    \node[anchor=north] (A) at (0,0) {\( V \)};
    \node[anchor=north] (B) at (-1,0) {\( V^* \)};
    \begin{knot}[draft mode=strands]
      \strand (A)
      .. controls +(0,\y) and +(0,\y) .. (B);
    \end{knot}
  \end{tikzpicture} \qquad
  \begin{tikzpicture}[xscale=0.75,baseline={(current bounding box.center)}]
    \node[anchor=north] (A) at (0,0) {\( V^* \)};
     \node[anchor=north] (B) at (-1,0) {\( V \)};
    \begin{knot}[draft mode=strands]
      \strand (A)
      .. controls +(0,-\y) and +(0,-\y) .. (B);
    \end{knot}
  \end{tikzpicture} \qquad
  \begin{tikzpicture}[xscale=0.75,baseline={(current bounding box.center)}]
    \node[anchor=north] (A) at (0,0) {\( V \)};
    \node[anchor=north] (B) at (1,0) {\( W \)};
    \node[anchor=south] (C) at (0,1) {\( W \)};
    \node[anchor=south] (D) at (1,1) {\( V \)};
    \begin{knot}[]
      \strand (B)
      .. controls +(-\x,\y) and +(\x,-\y) .. (C);
      \strand (A)
      .. controls +(\x,\y) and +(-\x,-\y) .. (D);
    \end{knot}
  \end{tikzpicture} \qquad
  \begin{tikzpicture}[xscale=1,baseline={(current bounding box.center)}]
    \node[anchor=north] (A) at (0,0) {\( {^*V} \)};
    \node[anchor=south] (B) at (0,2) {\( V^{*} \)};
    \draw (A)
    to [out=90, in=315] (-0.3,1)
    to [out=135, in=90, looseness=1.7] (-0.8,1);
    \draw[white,line width=4pt] ($(-0.3,1)+(-2pt,2pt)$) to ($(-0.3,1)+(2pt,-2pt)$);
    \draw (-0.8,1)
    to [out=270, in=225, looseness=1.7] (-0.3,1)
    to [out=45, in=270] (B);
  \end{tikzpicture} =
  \begin{tikzpicture}[xscale=0.75,baseline={(current bounding box.center)}]
    \node[anchor=north] (A) at (0,0) {\( {^*V} \)};
    \node[anchor=south] (B) at (0,2) {\( V^{*} \)};
    \draw (A) to coordinate[pos=0.5] (Circ) (B);
    \draw[fill=white, xscale=1.33] (Circ) circle (1.5pt);
  \end{tikzpicture}
  \caption{The evaluation, coevaluation, braiding morphism, and loop-de-loop depicted in the graphical calculus.  The loop-de-loop is abbreviated by a bead as indicated.}
  \label{fig-graphcalc}
\end{figure}

Let $\mathbf{1}$ denote the unit object in a rigid braided tensor category $\mathcal{A}$.  For any objects $V,W$ of a rigid braided tensor category, we have the following basic morphisms, of evaluation $\ev:V^*\ot V\to \mathbf{1}$, coevaluation $\coev:\mathbf{1}\to V\ot V^*$, and braiding $\sigma_{V,W}:V\ot W\to W\ot V$.  These are depicted in the graphical calculus in Figure \ref{fig-graphcalc}.

In addition, we will need the ``loop-de-loop'' morphism $l_V: {^*V}\to V^{*}$, defined as the composition,
$$l_V = (\ev_{^*V}\otimes \id_{V^*})\circ(\id_V \otimes\sigma_{V^*,{^*V}})\circ(\coev_V\otimes\id_{^*V}),$$ and depicted in the graphical calculus in Figure \ref{fig-graphcalc}.  To reduce clutter in diagrams, we will abbreviate this morphism by a bead, as indicated.

Following \cite{M95}, \cite{LM94}, \cite{M91}, the quantum coordinate algebra $\REAG$ arises as an algebra of (braided) matrix coefficients for the quantum group $U_q(\gl)$.  This equips $\REAG$, and hence $\REA$, with a number of structures we may exploit in our computations to come.  

Firstly, we recall that, for any finite dimensional (not necessarily irreducible) representation $W$ of $U_q(\gl)$, we have a \emph{canonical} homomorphism of $U_q(\gl)$-modules,
$$\iota_W: W^*\otimes W\to \REAG,$$
and if $W$ is a polynomial representation, then the image of $\iota_W$ lies in $\REA\subset \REAG$, recalling that $\REAG$ is obtained from $\REA$ by inverting the matrix coefficient of the quantum determinant representation.  Hence for a finite-dimensional polynomial representation $W$ of $U_q(\gl)$, let us denote by $C_{W}$ denote the image of $W^*\otimes W$ under the canonical $U_q(\gl)$-module homomorphism $\iota_W$.  

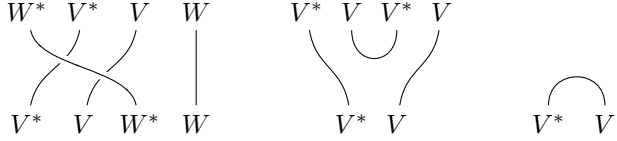
\begin{figure}[h]
  \centering
  \begin{tikzpicture}[xscale=0.75]
    \node[anchor=north] (A) at (0,0) {\( V^* \)};
    \node[anchor=north] (B) at (1,0) {\( V \)};
    \node[anchor=north] (C) at (2,0) {\( W^* \)};
    \node[anchor=north] (D) at (3,0) {\( W \)};

    \node[anchor=south] (Au) at (0,1) {\( W^* \)};
    \node[anchor=south] (Bu) at (1,1) {\( V^* \)};
    \node[anchor=south] (Cu) at (2,1) {\( V \)};
    \node[anchor=south] (Du) at (3,1) {\( W \)};
    \begin{knot}[]
      \flipcrossings{1,2}
      \strand (A)
      .. controls +(\x,\y) and +(-\x,-\y) .. (Bu);
      \strand (B)
      .. controls +(\x,\y) and +(-\x,-\y) .. (Cu);
      \strand (C)
      .. controls +(-\x,\y) and +(\x,-\y) .. (Au);
      \strand (D)
      .. controls +(0,\y) and +(0,-\y) .. (Du);
    \end{knot}
  \end{tikzpicture} \qquad
  \begin{tikzpicture}[xscale=0.6]
    \node[anchor=north] (A) at (1,0) {\( V^* \)};
    \node[anchor=north] (B) at (2,0) {\( V \)};

    \node[anchor=south] (Au) at (0,1) {\( V^* \)};
    \node[anchor=south] (Bu) at (1,1) {\( V \)};
    \node[anchor=south] (Cu) at (2,1) {\( V^* \)};
    \node[anchor=south] (Du) at (3,1) {\( V \)};
    \begin{knot}[]
      \strand (A)
      .. controls +(-\x,\y) and +(\x,-\y) .. (Au);
      \strand (B)
      .. controls +(\x,\y) and +(-\x,-\y) .. (Du);
      \strand (Bu)
      .. controls +(0,-\y) and +(0,-\y) .. (Cu);
    \end{knot}
  \end{tikzpicture} \qquad
  \begin{tikzpicture}[xscale=0.75]
    \node[anchor=north] (A) at (0,0) {\( V \)};
    \node[anchor=north] (B) at (-1,0) {\( V^* \)};
    \begin{knot}[draft mode=strands]
      \strand (A)
      .. controls +(0,\y) and +(0,\y) .. (B);
    \end{knot}
  \end{tikzpicture} 
\caption{The bi-algebra structure maps of multiplication, co-multiplication and co-unit in the $\REA$ are depicted in the graphical calculus.}\label{fig-bialg}
\end{figure}

The REA is a \emph{braided bi-algebra} in $\Rep_q(GL_N)$.  To specify the bialgebra structure maps on $\REA$, it suffices to specify them on each subspace $C_W$.  The multiplication $m$, co-multiplication $\Delta$, and co-unit $\epsilon$ morphisms are depicted in the graphical calculus, in Figure \ref{fig-bialg}.  The unit is just the identification $\mathbf{1}=C_{\mathbf{1}}$. We use the notational convention \( m^{(k)}:A^{\ot k} \longrightarrow A \) and \( \Delta^{(k)}:A\longrightarrow A^{\ot k} \) for the iterated product and coproduct. In particular \( m = m^{(2)} \) and \( \Delta = \Delta^{(2)} \).

\begin{figure}[h]
  \centering
    \[ A^2 = 
    \begin{tikzpicture}[xscale=0.6,
      baseline={([yshift=-0.5ex]current bounding box.center)}]

      \node[anchor=north,gray] at (-0.8,0.5) {\( \Delta \)};
      \node[anchor=south,gray] at (-0.8,1.5) {\(  m \)};
      
      \node[anchor=north] (A) at (1,0) {\( V^* \)};
      \node[anchor=north] (B) at (2,0) {\( V \)};

      \node[] (Au) at (0,1) {\( V^* \)};
      \node[] (Bu) at (1,1) {\( V \)};
      \node[] (Cu) at (2,1) {\( V^* \)};
      \node[] (Du) at (3,1) {\( V \)};

      \node[anchor=south] (Auu) at (0,2) {\( V^* \)};
      \node[anchor=south] (Buu) at (1,2) {\( V^* \)};
      \node[anchor=south] (Cuu) at (2,2) {\( V \)};
      \node[anchor=south] (Duu) at (3,2) {\( V \)};
      \begin{knot}[end tolerance=0em]
        \flipcrossings{1,2}
        \strand (A)
        .. controls +(-\x,\y) and +(\x,-\y) .. (Au);
        \strand (B)
        .. controls +(\x,\y) and +(-\x,-\y) .. (Du);
        \strand (Bu)
        .. controls +(0,-\y) and +(0,-\y) .. (Cu);

        \strand (Au)
        .. controls +(\x,\y) and +(-\x,-\y) .. (Buu);
        \strand (Bu)
        .. controls +(\x,\y) and +(-\x,-\y) .. (Cuu);
        \strand (Cu)
        .. controls +(-\x,\y) and +(\x,-\y) .. (Auu);
        \strand (Du)
        .. controls +(0,\y) and +(0,-\y) .. (Duu);
      \end{knot}
    \end{tikzpicture} =
    \begin{tikzpicture}[xscale=0.6,
      baseline={([yshift=-0.5ex]current bounding box.center)}]
      \node[anchor=north] (A) at (1,0) {\( V^* \)};
      \node[anchor=north] (B) at (2,0) {\( V \)};

      \node[anchor=south] (Au) at (0,1) {\( V^* \)};
      \node[anchor=south] (Bu) at (1,1) {\( V^* \)};
      \node[anchor=south] (Cu) at (2,1) {\( V \)};
      \node[anchor=south] (Du) at (3,1) {\( V \)}; 
      \begin{knot}[end tolerance=0em]
        \strand[] (Au)
        to [out=280, in=170] (0.5,0.7)
        to [out=350, in=260] coordinate[pos=0.4] (Circ) (Cu);
        \strand (B)
        .. controls +(\x,\y) and +(-\x,-\y) .. (Du);
        \strand (Bu)
        to [out=260, in=70] (0.5,0.7)
        to [out=250, in=100] (A);
      \end{knot}
      \draw[fill=white,xscale=1.666] (Circ) circle (1.5pt);
    \end{tikzpicture} \]
    \caption{The morphism $A^2:C_V\to \REA$ depicted here in the graphical calculus is a stand-in for the operation of squaring the matrix $A$.}\label{fig-A2}
\end{figure}
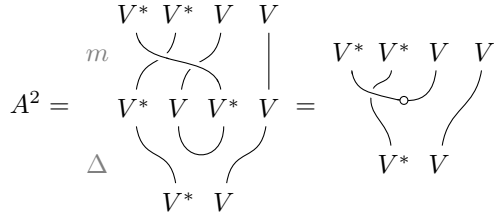

The matrix
$$ A = \sum_{i,j} E^i_j \ot a^j_i \in Mat_N(\CC) \otimes \REA \cong Mat_N(\REA)$$ defines an invariant element of \( \REA \ot V \ot V^* \), and thus may be identified with a homomorphism \(A\in \Hom_{U_q(\gl)}(C_V, \REA) \).  Note that we abuse notation and use $A$ for both its matrix and the associated homomorphism.  It is clear from the construction that
$$A=\iota_V: C_V \to\REA,$$
Likewise, $A^k$ corresponds both to the $k$th power of the matrix $A$, and to the morphism
$$A^k:C_V\xrightarrow{\iota_V} \REA \xrightarrow{\Delta^{(k)}} \REA^{\ot k} \xrightarrow{m^{(k)}} \REA,$$
where $\Delta^{(k)}$ and $m^{(k)}$ denote the $(k-1)$th iterated co-product, and product, respectively.  For example, the matrix,
$$A^2 = \sum_{i,j} E^i_j \ot \sum_k a^j_ka^k_i = \sum_{i,j} E^i_j \ot m \circ \Delta(a^i_j)$$ corresponds to the morphism depicted in Figure \ref{fig-A2}.

\begin{figure}[h]
  \[ A^k = 
    \begin{tikzpicture}[xscale=0.8,
      baseline={([yshift=-0.5ex]current bounding box.center)}]
      \node[anchor=north] (A) at (1,0) {\( V^* \)};
      \node[anchor=north] (B) at (2,0) {\( V \)};

      \node[anchor=south] (Au) at (0,2) {\( V^*  \)};
      \node[anchor=south] (Bu) at (1,2) {\( V^* \)};
      \node[anchor=south] (Cu) at (1.9,2) {\( V \)};
      \node[anchor=south] (Du) at (2.8,2) {\( V \)}; 
      \begin{knot}[end tolerance=0em, clip radius=8pt]
        \strand[double distance=12pt] ($(Au)+(0.25,-0.25)$)
        to [out=270, in=140] (0.5,1.5)
        to [out=320, in=180] (1.35,0.7)
        to [out=0, in=270] coordinate[pos=0] (Circ) ($(Cu)+(0.35,-0.25)$);
        
        \strand ($(Du) + (-0.0,-0.25)$)
        .. controls +(0,-1.7*\y) and +(0,+\y) .. (B);
        \strand ($(Bu)+(-0.2,-0.25)$)
        to [out=260, in=30] (0.3,1.5)
        to [out=210, in=135] (0.5,0.65)
        to [out=315, in=90] (A);
      \end{knot}
      \draw[double distance=4pt] ($(Au)+(0.25,-0.25)$)
        to [out=270, in=140] (0.5,1.5)
        to [out=320, in=180] (1.35,0.7)
        to [out=0, in=270] coordinate[pos=0] (Circ) ($(Cu)+(0.35,-0.25)$);
      \draw[fill=white,xscale=1.25] ($(Circ)+(0,-6pt)$) circle (1.5pt);
      \draw[fill=white,xscale=1.25] ($(Circ)+(0,-2pt)$) circle (1.5pt);
      \draw[fill=white,xscale=1.25] ($(Circ)+(0,2pt)$) circle (1.5pt);
      \draw[fill=white,xscale=1.25] ($(Circ)+(0,6pt)$) circle (1.5pt);

      \node[anchor=north] at ($(Cu)+(0.4,0)$) {\( \ldots \)};
      \node[anchor=north] at ($(Au)+(0.3,0)$) {\( \ldots \)};
      \begin{scope}[label box]
        \draw (-0.1,1.85) rectangle (1.1,1.2);
        \node[anchor=east] at (-0.1,1.525) {\( T_{(k\ldots 1)} \)};
      \end{scope}
    \end{tikzpicture} \]
  \caption{The formula for the morphism $A^k$ asserted in Lemma \ref{lem-Ak} is depicted here in the graphical calculus.}\label{fig-Ak}
\end{figure}
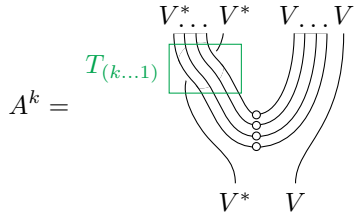

\begin{Lemma}\label{lem-Ak} The morphism $A^k: C_V \to C_{V^{\ot k}}\subseteq \REA$ is as depicted in Figure \ref{fig-Ak}.
\end{Lemma}

\begin{proof} The proof of this is an easy induction starting from $k=2$, with the induction step illustrated in Figure \ref{fig-Akpf}.

  \begin{figure}[h]
    \[ A^{k+1} = 
    \begin{tikzpicture}[xscale=0.8,
      baseline={([yshift=-0.5ex]current bounding box.center)}]
      \node[anchor=north] (Ab) at (1,-1) {\( V^* \)};
      \node[anchor=north] (Bb) at (2,-1) {\( V \)};

      \coordinate[anchor=north] (A) at (1,0);
      \coordinate[anchor=north] (B) at (2,0);

      \node[anchor=south] (Au) at (0,2) {};
      \node[anchor=south] (Bu) at (1,2) {};
      \node[anchor=south] (Cu) at (1.9,2) {};
      \node[anchor=south] (Du) at (2.8,2) {};

      \node[anchor=south] (Auu) at (0,2.5) {\( V^*  \)};
      \node[anchor=south] (Buu) at (1,2.5) {\( V^* \)};
      \node[anchor=south] (Cuu) at (1.9,2.5) {\( V \)};
      \node[anchor=south] (Duu) at (2.8,2.5) {\( V \)};

      \node[anchor=south] (eAuu) at (-0.6,2.5) {\( V^* \)};
      \node[anchor=south] (eDuu) at (3.3,2.5) {\( V \)};
      \begin{knot}[end tolerance=0em, clip radius=8pt]
        \strand[double distance=12pt] ($(Auu)+(0.25,-0.25)$)
        to [out=270, in=90] ($(Au)+(0.25,-0.25)$)
        to [out=270, in=140] (0.5,1.5)
        to [out=320, in=180] (1.35,0.7)
        to [out=0, in=270] coordinate[pos=0] (Circ) ($(Cu)+(0.35,-0.25)$)
        to [out=90, in=270] ($(Cuu)+(0.35,-0.25)$);
        
        \strand ($(Duu) + (-0.0,-0.25)$)
        to [out=270, in=90] ($(Du) + (-0.0,-0.25)$)
        .. controls +(0,-1.7*\y) and +(0,+\y) .. (B);
        \strand ($(Buu) + (-0.2,-0.25)$)
        to [out=270, in=90] ($(Bu)+(-0.2,-0.25)$)
        to [out=260, in=30] (0.3,1.5)
        to [out=210, in=135] (0.5,0.65)
        to [out=315, in=90] (A) to (Ab);

        \strand (eDuu)
        to [out=270, in=90] ($(eDuu)+(0,-2.6)$)
        to [out=270, in=90, looseness=0.9] (Bb);
      \end{knot}
      \draw[double distance=4pt] ($(Auu)+(0.25,-0.25)$)
      to [out=270, in=90] ($(Au)+(0.25,-0.25)$)
      to [out=270, in=140] (0.5,1.5)
      to [out=320, in=180] (1.35,0.7)
      to [out=0, in=270] coordinate[pos=0] (Circ) ($(Cu)+(0.35,-0.25)$)
      to [out=90, in=270] ($(Cuu)+(0.35,-0.25)$);
      
      \draw[fill=white,xscale=1.25] ($(Circ)+(0,-6pt)$) circle (1.5pt);
      \draw[fill=white,xscale=1.25] ($(Circ)+(0,-2pt)$) circle (1.5pt);
      \draw[fill=white,xscale=1.25] ($(Circ)+(0,2pt)$) circle (1.5pt);
      \draw[fill=white,xscale=1.25] ($(Circ)+(0,6pt)$) circle (1.5pt);

      \draw[line width=3.5pt, white] (B)
      to [out=270, in=270, looseness=1.5] ($(eDuu)+(-0.2,-2.3)$)
      to ($(eDuu)+(-0.2,-0.9)$)
      to [out=90, in=270, looseness=0.4] (eAuu);
      
      \draw (B)
      to [out=270, in=270, looseness=1.5] ($(eDuu)+(-0.2,-2.3)$)
      to ($(eDuu)+(-0.2,-0.9)$)
      to [out=90, in=270, looseness=0.4] (eAuu);

      \node[anchor=north] at ($(Cuu)+(0.4,0)$) {\( \ldots \)};
      \node[anchor=north] at ($(Auu)+(0.3,0)$) {\( \ldots \)};

      \begin{scope}[label box]
        \draw (-0.7,2.4) rectangle (3.45,1.95);
        \node[anchor=east] at (-0.7,2.165) {\( m \)};
        
        \draw (-0.15,1.85) rectangle (2.95,0.2);
        \node[anchor=east] at (-0.15,1.025) {\( A^k \)};
        
        \draw (0.8,0.1) rectangle (3.45,-0.8);
        \node[anchor=east] at (0.8,-0.35) {\( \Delta \)};
      \end{scope}
    \end{tikzpicture} =
    \begin{tikzpicture}[xscale=0.8,
      baseline={([yshift=-0.5ex]current bounding box.center)}]
      \node[anchor=north] (A) at (1,0) {\( V^* \)};
      \node[anchor=north] (B) at (2,0) {\( V \)};

      \node[anchor=south] (Au) at (0,2) {\( V^*  \)};
      \node[anchor=south] (Bu) at (1,2) {\( V^* \)};
      \node[anchor=south] (Cu) at (1.9,2) {\( V \)};
      \node[anchor=south] (Du) at (2.8,2) {\( V \)}; 
      \begin{knot}[end tolerance=0em, clip radius=8pt]
        \strand[double distance=12pt] ($(Au)+(0.25,-0.25)$)
        to [out=270, in=140] (0.5,1.5)
        to [out=320, in=180] (1.35,0.7)
        to [out=0, in=270] coordinate[pos=0] (Circ) ($(Cu)+(0.35,-0.25)$);
        
        \strand ($(Du) + (-0.0,-0.25)$)
        .. controls +(0,-1.7*\y) and +(0,+\y) .. (B);
        \strand ($(Bu)+(-0.2,-0.25)$)
        to [out=260, in=30] (0.3,1.5)
        to [out=210, in=135] (0.5,0.65)
        to [out=315, in=90] (A);
      \end{knot}
      \draw[double distance=4pt] ($(Au)+(0.25,-0.25)$)
      to [out=270, in=140] (0.5,1.5)
      to [out=320, in=180] (1.35,0.7)
      to [out=0, in=270] coordinate[pos=0] (Circ) ($(Cu)+(0.35,-0.25)$);
      \draw[fill=white,xscale=1.25] ($(Circ)+(0,-6pt)$) circle (1.5pt);
      \draw[fill=white,xscale=1.25] ($(Circ)+(0,-2pt)$) circle (1.5pt);
      \draw[fill=white,xscale=1.25] ($(Circ)+(0,2pt)$) circle (1.5pt);
      \draw[fill=white,xscale=1.25] ($(Circ)+(0,6pt)$) circle (1.5pt);

      \node[anchor=north] at ($(Cu)+(0.4,0)$) {\( \ldots \)};
      \node[anchor=north] at ($(Au)+(0.3,0)$) {\( \ldots \)};

      \begin{scope}[label box]
        \draw (-0.1,1.85) rectangle (1.1,1.2);
        \node[anchor=east] at (-0.1,1.525) {\( T_{(k+1\ldots 1)} \)};
      \end{scope}
      
    \end{tikzpicture} \]
    \caption{The morphism $A^{k+1}$ is by definition the composition of $m\circ A^k\circ\Delta$.  Inductively assuming our formula for $A^k$, we obtain the formula for $A^{k+1}$ by elementary isotopies of braids.}\label{fig-Akpf}
  \end{figure}
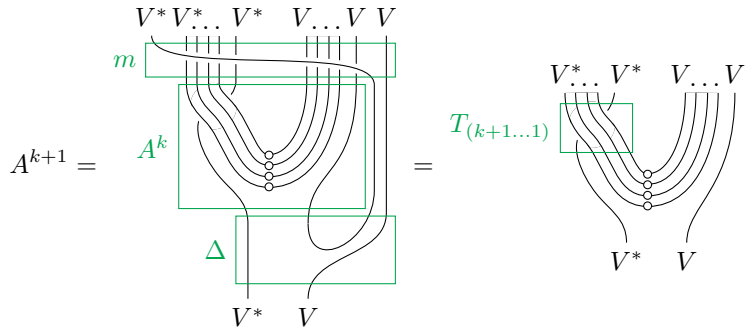
\end{proof}

\subsection{The finite Hecke algebra}
\label{sec:finite-hecke-algebra}
In order to formulate Reshetikhin's central elements, we need to recall some basic facts about the finite Hecke algebra.

\begin{Definition}
The finite Hecke algebra, $H_q(n)$ is generated by the symbols
$$T_1,\ldots, T_{n-1},$$
subject to the relations
$$T_iT_j=T_jT_i \textrm{ if $|i-j|\geq 2$}, \qquad T_iT_{i+1}T_i = T_{i+1}T_iT_{i+1} \textrm{ for $i=1,\ldots n-1$}$$
$$ T_i-T_{i}^{-1} = q-q^{-1}, \textrm{ for $i=1,\ldots, n-1$}.$$
\end{Definition}

\begin{Definition} Given an element $\sigma\in S_n$, we denote by $T_\sigma$ its canonical lift to $H_q(n)$, obtained by writing $\sigma=s_{i_1}\ldots s_{i_\ell}$ as a reduced word in the generators $s_i$ and lifting those to $T_\sigma = T_{i_1}\ldots T_{i_\ell}$.
\end{Definition}
The element $T_\sigma$ depends only on $\sigma$, not on the reduced word presentation of $\sigma$ chosen.  In particular, given a cycle of the form $(k\ldots i)$, for $i<k$ we denote by
$$T_{(k\ldots i)} = T_{i}\cdots T_{k-1}$$
the lift obtained in this way.  The braid $T_{(k\ldots i)}$ is depicted in Figure \ref{fig-Ak}.

\begin{Definition}
We define elements $\omega_k, \overline{\omega}_k\in H_q(n)$ as:
$$\omega_k = \sum_{\sigma\in S_k} (-q)^{-\ell(\sigma)} T_\sigma \in H_q(k),\qquad \overline{\omega}_k = \frac{\omega_k}{[k]_q!}.$$
where we denote the quantum integers, and the quantum factorial, by:
$$[k]_q=1+q^{-2} + \ldots + q^{-2(k-1)}=\frac{1-q^{-2k}}{1-q^{-2}}, \qquad [k]_q! = [k]_q[k-1]_q\cdots [2]_q[1]_q.$$
\end{Definition}

For $k<n$, we abuse notation and regard $\omega_k,\overline{\omega}_k$ as an element of $H_q(n)$ via the obvious inclusion $H_q(k)\subset H_q(n)$, sending $T_i$ to $T_i$.

\begin{Remark} We note that $\omega_k$ is only a weak idempotent, as $\omega_k\omega_k = [k]_q! \omega_k$, while its normalization $\overline{\omega}_k$ is an idempotent.\end{Remark}
 
Recall that, for $k<N$, quantum Schur-Weyl duality \cite[Section~8.6.3]{KS97} gives an isomorphism
$$\End(V^{\ot k})\cong H_q(k),$$
sending the generators $T_i$ to the braiding on the $i$th and $(i+1)$st components.  In particular, the fundamental representations of $U_q(\mathfrak{gl}_N)$ arise as sub-modules of $V^{\ot k}$ obtained by projecting along the weak idempotents $\omega_k$.  We use this formulation in the next section.

\subsection{Reshetikhin's central elements}
\label{sec:resh-centr-elem}
Finally, let us turn to the construction of Reshetikhin's central elements $c_k$.  For any representation $W$ in $\Rep_q(GL_N)$, we define $\tr_q(W) = \iota_W\circ\sigma_{W,W^*}\circ\coev_W(1)$.  This is depicted in Figure \ref{fig-trqA}.

\begin{figure}[h]
  \centering
  \begin{tikzpicture}[xscale=0.75,baseline={([yshift=-0.5ex]current bounding box.center)}]
    \node[anchor=south] (A) at (0,0) {\( W \)};
    \node[anchor=south] (B) at (-1,0) {\( W^* \)};
    \coordinate (xs) at (-0.5,-0.5);
    \coordinate (min) at (-0.5,-0.8);
    \draw (A)
    to [out=270, in=45] (xs)
    to [out=225, in=180, looseness=1.7] (min);
    \draw[white,line width=4pt] ($(xs)+(-2pt,2pt)$) to ($(xs)+(2pt,-2pt)$);
    \draw (min)
    to [out=0, in=315, looseness=1.7] (xs)
    to [out=135, in=270] (B);
  \end{tikzpicture} =
  \begin{tikzpicture}[xscale=0.75,baseline={([yshift=-0.5ex]current bounding box.center)}]
    \node[anchor=south] (A) at (0,0) {\( W \)};
    \node[anchor=south] (B) at (-1,0) {\( W^* \)};
    \begin{knot}[draft mode=strands]
      \strand (A)
      .. controls +(0,-\y) and +(0,-\y) .. coordinate[pos=0.5] (Circ) (B);
    \end{knot}
    \draw[fill=white, xscale=1.33] (Circ) circle (1.5pt);
  \end{tikzpicture}
  \caption{The morphism $tr_q(W)$ defines a canonical invariant in $\REAG$, for every $W\in\Rep_q(GL_N)$.}\label{fig-trqA}
\end{figure}

\begin{Definition}
Reshetikhin's canonical central elements are:
$$c_k = \tr_q(\Lambda_q^k(V)) = \tr_q(L(\omega_k)),$$
i.e. the central invariant corresponding to the irreducible representation of highest weight $\omega_k$, equivalently the $k$th exterior power of the defining representation.
\end{Definition}

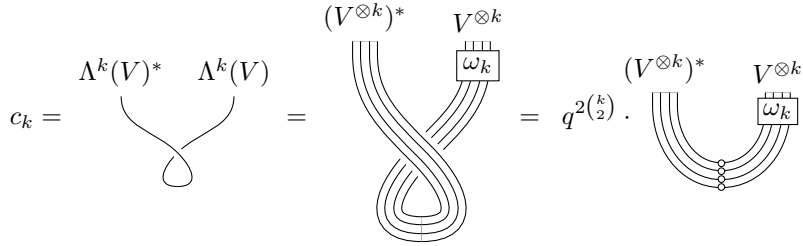
\begin{figure}[h]
  \centering \[ c_k =   
  \begin{tikzpicture}[xscale=1,yscale=1,baseline={([yshift=-0.5ex]current bounding box.center)}]
    \node[anchor=south] (A) at (0,0) {\( \Lambda^k (V) \)};
    \node[anchor=south] (B) at (-1.5,0) {\( \Lambda^k (V)^* \)};
    \coordinate (xs) at (-0.75,-0.75);
    \coordinate (min) at (-0.75,-1.2);
    \draw (A)
    to [out=270, in=45] (xs)
    to [out=225, in=180, looseness=1.7] (min);
    \draw[white,line width=4pt] ($(xs)+(-3pt,3pt)$) to ($(xs)+(3pt,-3pt)$);
    \draw (min)
    to [out=0, in=315, looseness=1.7] (xs)
    to [out=135, in=270] (B);
  \end{tikzpicture} =
  \begin{tikzpicture}[xscale=1,yscale=1,baseline={([yshift=-0.5ex]current bounding box.center)}]
    \node[anchor=south] (A) at (0,0.5) {\( V^{\otimes k} \)};
    \node[anchor=south] (B) at (-1.5,0.5) {\( (V^{\otimes k})^* \)};
    \coordinate (xs) at (-0.75,-1);
    \coordinate (min) at (-0.75,-2);
    \draw[double distance=9pt] (A) to ($(A)+(0,-0.5)$)
    to [out=270, in=45] (xs)
    to [out=225, in=180, looseness=1.7] (min);
    \draw[double distance=3pt] (A) to ($(A)+(0,-0.5)$)
    to [out=270, in=45] (xs)
    to [out=225, in=180, looseness=1.7] (min);
    \draw[white,line width=14pt] ($(xs)+(-20pt,20pt)$) to ($(xs)+(20pt,-20pt)$);
    \draw[double distance=9pt] (min)
    to [out=0, in=315, looseness=1.7] (xs)
    to [out=135, in=270] ($(B)+(0,-0.5)$)
    to (B);
    \draw[double distance=3pt] (min)
    to [out=0, in=315, looseness=1.7] (xs)
    to [out=135, in=270] ($(B)+(0,-0.5)$)
    to (B);
    \draw[fill=white] ($(A)+(-8pt,-0.4)$) rectangle ($(A)+(8pt,-0.8)$);
    \node at ($(A)+(0,-0.6)$) {\( \omega_k \)};
  \end{tikzpicture} = \;\; q^{2\binom{k}{2}} \cdot \!\!\!\!\!
  \begin{tikzpicture}[xscale=1,yscale=1,baseline={([yshift=-0.5ex]current bounding box.center)}]
    \node[anchor=south] (A) at (0,0.5) {\( V^{\otimes k} \)};
    \node[anchor=south] (B) at (-1.5,0.5) {\( (V^{\otimes k})^* \)};
    \coordinate (min) at (-0.75,-0.6);
    \draw[double distance=9pt] (A) to ($(A)+(0,-0.5)$)
    to [out=270, in=0, looseness=1.0] (min)
    to [out=180, in=270, looseness=1.0] ($(B)+(0,-0.5)$)
    to (B);
    \draw[double distance=3pt] (A) to ($(A)+(0,-0.5)$)
    to [out=270, in=0, looseness=1.0] (min)
    to [out=180, in=270, looseness=1.0] ($(B)+(0,-0.5)$)
    to (B);
    \draw[fill=white] ($(A)+(-7.5pt,-0.35)$) rectangle ($(A)+(7.5pt,-0.7)$);
    \node at ($(A)+(0,-0.525)$) {\( \omega_k \)};

    \draw[fill=white] ($(min)+(0,-4.5pt)$) circle (1.25pt);
    \draw[fill=white] ($(min)+(0,-1.5pt)$) circle (1.25pt);
    \draw[fill=white] ($(min)+(0,1.5pt)$) circle (1.25pt);
    \draw[fill=white] ($(min)+(0,4.5pt)$) circle (1.25pt);
  \end{tikzpicture} \]
\caption{The canonical central elements are the quantum traces of the $k$th exterior power of the defining representation, which may be embedded into $C_{V^{\ot k}}$ via the corresponding projector $\omega_k$.  See Lemma \ref{lem-ckid}.}\label{fig-ckid}
\end{figure}

\begin{Lemma}\label{lem-ckid} We have the identities for Reshetikhin's elements depicted in Figure~\ref{fig-ckid}.\end{Lemma}
\begin{proof}
The first equality is the definition of $c_k$.  The second equality follows from the construction of $\Lambda^k(V)$ as a submodule of $V^{\ot k}$, determined as the image of the projector $\omega_k$ in the finite Hecke algebra, which under Schur-Weyl duality identifies with endomorphisms of $V^{\ot k}$.  The third equality follows by decomposing the $k$-fold loop-de-loop into a k-fold inverse double-braiding, composed with the individual loop-de-loops on each strand.  The $k$-fold inverse double-braiding has $2\cdot {k \choose 2}$ inverse crossings, each of which simplifies to $-q$ upon meeting the projector $\pi$.
\end{proof}

\begin{Theorem}[{\cite[Theorem 2]{B98}}]
The elements $c_1,\ldots, c_N$ generate the center of $\REA$, and together with the inverse of $c_N$ they generate the center of $\REAG$.
\end{Theorem}

\begin{Remark}
We note in passing that the appearance of weight functions $\wt(I)$ in our formulas can ultimately be traced back to the appearance of loop-de-loop's appearing in the right-hand side of Figure \ref{fig-ckid}.  It is well-known that the loop-de-loop on the defining representation acts diagonally on the weight basis with powers $1,q^{-2},\cdots, q^{2-N}$.
\end{Remark}

\section{Quantum Girard-Newton identities}
\label{sec:quant-girard-newton}
In this section, we prove a quantum Newton identity relating the canonical generators $c_k$ to the quantum trace generators $s_k$.  As a corollary, we identify Reshetkhin's canonical generators $c_k$ with (a scalar multiple of) the generators $\sigma_k$ of \cite{PS96,GPS97}.  As a corollary, we can apply the quantum Cayley-Hamilton identity from \cite{GPS97} to the $c_k$'s, as asserted in Theorem~\ref{main-thm}.

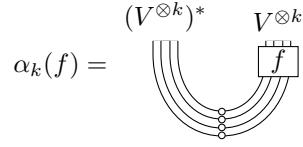
\begin{figure}[h]
  \centering \[ \alpha_k(f) =   
  \begin{tikzpicture}[xscale=1,yscale=1,baseline={([yshift=-0.5ex]current bounding box.center)}]
    \node[anchor=south] (A) at (0,0.5) {\( V^{\otimes k} \)};
    \node[anchor=south] (B) at (-1.5,0.5) {\( (V^{\otimes k})^* \)};
    \coordinate (min) at (-0.75,-0.6);
    \draw[double distance=9pt] (A) to ($(A)+(0,-0.5)$)
    to [out=270, in=0, looseness=1.0] (min)
    to [out=180, in=270, looseness=1.0] ($(B)+(0,-0.5)$)
    to (B);
    \draw[double distance=3pt] (A) to ($(A)+(0,-0.5)$)
    to [out=270, in=0, looseness=1.0] (min)
    to [out=180, in=270, looseness=1.0] ($(B)+(0,-0.5)$)
    to (B);
    \draw[fill=white] ($(A)+(-7.5pt,-0.35)$) rectangle ($(A)+(7.5pt,-0.75)$);
    \node at ($(A)+(0,-0.55)$) {\( f \)};

    \draw[fill=white] ($(min)+(0,-4.5pt)$) circle (1.25pt);
    \draw[fill=white] ($(min)+(0,-1.5pt)$) circle (1.25pt);
    \draw[fill=white] ($(min)+(0,1.5pt)$) circle (1.25pt);
    \draw[fill=white] ($(min)+(0,4.5pt)$) circle (1.25pt);
  \end{tikzpicture} \]
  \caption{The linear map $\alpha_k(f)$ of Definition \ref{def:alphak} is depicted in the graphical calulus.}
  \label{fig:alphak}
\end{figure}

\begin{Definition}\label{def:alphak} We define a linear map,
$$\alpha_k:  \End_{U_q(\gl)}(V^{\ot k}) \to \REA,$$
by sending $f\in \End(V^{\ot k})$ to the composition,
$$\alpha_k(f) = (l_{V}^{\ot k}\ot f)\circ\coev_{(^*V)^{\ot k}},$$
depicted in Figure \ref{fig:alphak}.
\end{Definition}

\begin{Proposition} The representation,
\begin{align*}\rho: H_q(n) &\to \End(V^{\otimes n}),\\
T_i &\mapsto (\sigma_{V,V})_{i,i+1},\end{align*}
descends through $\alpha_{n}$ to a linear map,
$$\HH(H_q(n))\to C_{V^{\otimes n}}.$$
\end{Proposition}

\begin{proof}
Amongst the defining linear relations of $C_{V^{\otimes n}}$ are the relations 
$$\phi^*f\otimes v = f\otimes\phi(v),$$
for any $\phi:V^{\otimes n}\to V^{\otimes n}$, in particular for the braiding morphisms.  Under $\alpha_{n}$, these identify the operation of post-composing and pre-composing with the map $\phi$. Hence $\alpha_{n}\circ \rho$ respects the defining linear relations, $ab-ba$, of $\HH$.
\end{proof}

\begin{Proposition}\label{QNI-prop}
We have the following identity in $\HH(H_q(k))$:
\begin{align*}\omega_k &= \omega_{k-1} + (-q^{-1})[k-1]_qT_{k-1}\omega_{k-1}.
\end{align*}

\end{Proposition}
\begin{proof}
By definition, we have:
$$\omega_k = \omega_{k-1} + \sum_{i=1}^{k-1} (-q)^{k-i}T_{(k\ldots i)} \omega_{k-1}.$$
However, in $\HH$, we can pull each $T_j$, with $j<k-1$, appearing in the product
$T_{(k\cdots i)} = T_i\cdots T_{k-1}$
around to the right, where it multiplies against $\omega_{k-1}$ to give a further factor of $(-q^{-1})$.  The claimed formula follows immediately.
\end{proof}

\begin{Corollary}\label{cor-HeckeQNI} The following identity holds in $H_q(n)$:

\begin{align*}[k]_q\overline{\omega}_k 
&= \sum_{j=1}^{k} (-q^{-1})^{k-j-1}\overline{\omega}_{j-1}T_{(j\cdots k)}
\end{align*}
\end{Corollary}

\begin{proof} This is a simple induction starting from Proposition \ref{QNI-prop}.\end{proof}

\begin{Theorem}
  \label{thm:qNewton-id}
  We have the following ``quantum Newton identity'':
\begin{equation}
  [k]_q s_k = \sum_{j=1}^{k} \frac{c_{j-1}}{[j-1]_q!}s_{k-j}.
\end{equation}
\end{Theorem}

\begin{proof}
The theorem follows from Corollary \ref{cor-HeckeQNI}, as soon as we establish the following identity:
\begin{equation}\alpha_{k} (\omega_{j-1}T_{(k\ldots j)}) = \alpha_{k}(\omega_{j-1})\alpha_{k}(T_{(k\ldots j)}) = c_{j-1}s_{k-j}.\label{eqn-alphaidentity}\end{equation}

First, we recall the identity for $A^k$ given in Lemma \ref{lem-Ak}.  Note that pre-composing this identity with with $\tr_q: \mathbf{1} \to C_V$ immediately yields:
$$s_k=\tr_q(A^k) = \alpha_{k}(T_{(k\ldots 1)}).$$

On the other hand, Reshetikhin's central elements $c_k$ were defined the as the quantum trace of the $k$th fundamental representation, hence we have:
$$c_k = \tr_q(\Lambda^k(V)) = \alpha_{k}(\omega_k)  .$$

While $\alpha_{k}$ is not an algebra homomorphism, we note that when two braids $a$ and $b$ involve disjoint strands $1,\ldots, i$ and $i+1,\ldots i+j$, then we have that $\alpha_{i+j}(ab) = \alpha_i(a)\alpha_j(b)$.  This is depicted in Figure \ref{fig-ctimesp}.

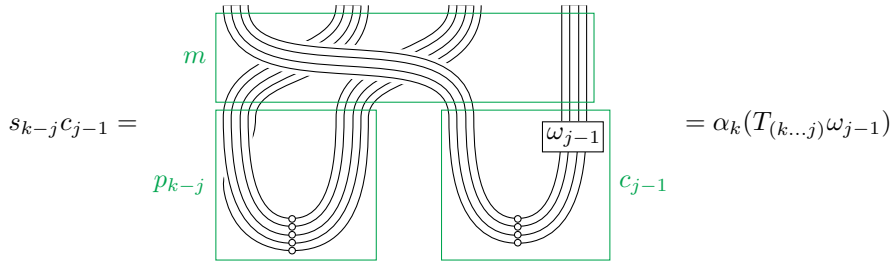
\begin{figure}[h]
  \centering
  \[ s_{k-j}c_{j-1} =   
    \begin{tikzpicture}[xscale=1,yscale=1,baseline={([yshift=-0.5ex]current bounding box.center)}]
      \def\hspc{1.5}
      \def\vspc{1.5}
      \def\sepr{7.5pt}
      \coordinate (A) at (0,0);
      \coordinate (B) at ($(A)+(\hspc,0)$);
      \coordinate (C) at ($(B)+(\hspc,0)$);
      \coordinate (D) at ($(C)+(\hspc,0)$);
      
      \coordinate (min1) at ($0.5*(A)+0.5*(B)+(0,-2*\vspc)$);
      \coordinate (min2) at ($0.5*(C)+0.5*(D)+(0,-2*\vspc)$);

      \coordinate (exit) at (-0.15,-2.5);
      \coordinate (dirc) at (-0.2,1);

      \draw ($(B)+(\sepr,0)$)
      .. controls +(0,-1.05) and +(0,0.6) .. ($(A)+(0,-\vspc)+(\sepr,0)$)
      .. controls +(0,-0.6) and +($0.5*(dirc)$) .. (exit)
      .. controls +($-0.75*(dirc)$) and +(-0.1,0) .. ($(min1)+(0,-\sepr)$)
      .. controls +(1,0) and +(0,-0.6) .. ($(B)+(0,-\vspc)+(\sepr,0)$)
      .. controls +(0,0.6) and +(0,-1.05) .. ($(C)+(\sepr,0)$);

      \draw[white,line width=12.5pt] (B)
      to [out=270, in=90] ($(A)+(0,-\vspc)$)
      to [out=270, in=180] (min1)
      to [out=0, in=270] ($(B)+(0,-\vspc)$)
      to [out=90, in=270] (C);
      
      \draw[double distance=9pt] (B)
      to [out=270, in=90] ($(A)+(0,-\vspc)$)
      to [out=270, in=180] (min1)
      to [out=0, in=270] ($(B)+(0,-\vspc)$)
      to [out=90, in=270] (C);

      \draw[double distance=3pt] (B)
      to [out=270, in=90] ($(A)+(0,-\vspc)$)
      to [out=270, in=180] (min1)
      to [out=0, in=270] ($(B)+(0,-\vspc)$)
      to [out=90, in=270] (C);

      \draw[white,line width=13pt] (A)
      to[out=270, in=90] ($(C)+(0,-\vspc)$)
      to[out=270, in=180] (min2)
      to[out=0, in=270] ($(D)+(0,-\vspc)$)
      to[out=90, in=270] (D);

      \draw[double distance=9pt] (A)
      to[out=270, in=90] ($(C)+(0,-\vspc)$)
      to[out=270, in=180] (min2)
      to[out=0, in=270] ($(D)+(0,-\vspc)$)
      to[out=90, in=270] (D);

      \draw[double distance=3pt] (A)
      to[out=270, in=90] ($(C)+(0,-\vspc)$)
      to[out=270, in=180] (min2)
      to[out=0, in=270] ($(D)+(0,-\vspc)$)
      to[out=90, in=270] (D);

      \coordinate (rect) at ($(D)+(0,-\vspc)+(0,-0.25)$);

      \draw[fill=white] ($(rect)+(-12pt,0.2)$) rectangle ($(rect)+(11pt,-0.2)$);
      \node at (rect) {\( \omega_{j-1} \)};

      \draw[fill=white] ($(min1)+(0,-4.5pt)$) circle (1.25pt);
      \draw[fill=white] ($(min1)+(0,-1.5pt)$) circle (1.25pt);
      \draw[fill=white] ($(min1)+(0,1.5pt)$) circle (1.25pt);
      \draw[fill=white] ($(min1)+(0,4.5pt)$) circle (1.25pt);
      \draw[fill=white] ($(min1)+(0,-7.5pt)$) circle (1.25pt);

      \draw[fill=white] ($(min2)+(0,-4.5pt)$) circle (1.25pt);
      \draw[fill=white] ($(min2)+(0,-1.5pt)$) circle (1.25pt);
      \draw[fill=white] ($(min2)+(0,1.5pt)$) circle (1.25pt);
      \draw[fill=white] ($(min2)+(0,4.5pt)$) circle (1.25pt);

      \begin{scope}[label box]
        \coordinate (mrecttl) at ($(A)+(-7.5pt,-3pt)$);
        \coordinate (mrectbr) at ($(D)+(7.5pt,-\vspc)+(0,6pt)$);
        \draw (mrecttl) rectangle (mrectbr);
        \gettikzxy{(mrectbr)}{\mrectbrx}{\mrectbry}
        \gettikzxy{(mrecttl)}{\mrecttlx}{\mrecttly}
        \node[anchor=east] at ($0.5*(mrecttl)+0.5*(\mrecttlx,\mrectbry)$) {\( m \)};
        
        \coordinate (precttl) at ($(A)+(0,-\vspc)+(-7.5pt,+3pt)$);
        \coordinate (prectbr) at ($(B)+(10.5pt,-2.4*\vspc)+(0,6pt)$);
        \draw (precttl) rectangle (prectbr);
        \gettikzxy{(prectbr)}{\prectbrx}{\prectbry}
        \gettikzxy{(precttl)}{\precttlx}{\precttly}
        \node[anchor=east] at ($0.5*(precttl)+0.5*(\precttlx,\prectbry)$) {\( p_{k-j} \)};

        \coordinate (crecttl) at ($(C)+(0,-\vspc)+(-7.5pt,+3pt)$);
        \coordinate (crectbr) at ($(D)+(13.5pt,-2.4*\vspc)+(0,6pt)$);
        \draw (crecttl) rectangle (crectbr);
        \gettikzxy{(crectbr)}{\crectbrx}{\crectbry}
        \gettikzxy{(crecttl)}{\crecttlx}{\crecttly}
        \node[anchor=west] at ($0.5*(crectbr)+0.5*(\crectbrx,\crecttly)$) {\( c_{j-1} \)};
      \end{scope}
      
    \end{tikzpicture} = \alpha_k(T_{(k\ldots j)}\omega_{j-1}) \]
  \caption{Because $T_{(k\cdots j)}$ and $\omega_{j-1}$ involve disjoint strands, we have the identity $\alpha_k(T_{(k\cdots j)}\omega_{j-1}) = \alpha_{k-j}(T_{(k\cdots j)})\alpha_{j-1}(\omega_{j-1})$}\label{fig-ctimesp}
\end{figure}

Hence Equation \eqref{eqn-alphaidentity} follows, and the theorem is proved.
\end{proof}

As a consequence, we have that Reshetikhin's central generators $c_k$ coincide with Gurevich, Pyatov and Saponov's generators $\sigma_k$ (up to an inconsequential scaling), since both sets of generators are uniquely determined by the fact that they satisfy a quantum Newton identity with respect to quantum power traces, as in the theorem.  Since the quantum Cayley-Hamilton theorem was proved already in \cite{GPS97} for the $\sigma_k$'s, we immediately obtain the analogous quantum Cayley-Hamilton identity asserted in the introduction, for the $c_k$'s.

\begin{Corollary}
  \label{cor:qCH}
The central elements \( c_k \) satisfy the the following quantum Cayley Hamilton identity
\begin{equation}
  \sum_{k=0}^N(-q^2)^{N-k}c_{N-k}\cdot A^k =0.
\end{equation}
\end{Corollary}

\begin{Remark}
It is possible to use similar techniques as above to directly obtain the quantum Cayley-Hamilton formula involving Reshetikhin's generators.  However, as it is already proved in \cite{GPS97} that the QCH identities follows from the quantum Girard-Newton identities, we omit such a discussion here.
\end{Remark}

\section{Quantum minors for the FRT and the REA}
\label{sec:minors-combinatorics}

In this section we define three of our primary objects of study: the quantum minors for the FRT and the REA algebra, and the truncated minors in the REA.  We recall a well-known row expansion formula for quantum minors in the FRT, and we formulate an analogous one for the truncated quantum minors.

\subsection{Permutation statistics}
\label{sec:permutations}

First we recall the classical definitions for permutations in \( \SN \), that is, bijections from \( [N] \) to itself. The \emph{length} \( \lng(\sigma) \) of a permutation \( \sigma \in \SN \) is the number of inversions, that is
\[ \lng(\sigma) = \# \setc{(i,j) \in [N]^2}{i < j \text{ and } \sigma(i) > \sigma(j)}. \]
We also recall the definition of the \emph{exceedance} \( e(\sigma) \) of a permutation \( \sigma \in \SN \) as
\[ e(\sigma) = \# \setc{i \in [N]}{\sigma(i) > i}. \]

We extend the definitions of the above statistics to bijections to more general bijections. We fix subsets \( I,J \subset [N] \) such that \( \# I = \# J \), and an \emph{auxiliary set} \( U \subset (I\cup J)^\comp \) (the complement of \( I\cup J \) in \( [N] \)). Given a bijection \( \tau\map{I}{J} \) let \( \hat{\tau}\map{I\cup U}{J \cup U} \) be the bijection defined by
\begin{equation*}
  \hat{\tau}(i) =
  \begin{cases}
    \tau(i) &\text{if } i \in I, \\
    i       &\text{otherwise}.
  \end{cases}
\end{equation*}

\begin{Definition}
  \label{def:length-exceedance}
The \emph{exceedance}, \( e(\tau) \) and the \emph{length}, \( \lng_U(\tau) \), of a bijection \( \tau\map{I}{J} \) are defined by
\begin{align*}
  e(\tau) &= \# \setc{a \in I}{\tau(a) > a}, \text{ and } \\
  \lng_U(\tau) &= \# \setc{(a,b) \in (I\cup U) \times (I \cup U)}{a<b \text{ and } \hat{\tau}(a) > \hat{\tau}(b)}.
\end{align*}
\end{Definition}

Note that the exceedance does not depend on the auxiliary set \( U \). In the case \( I = J \) and \( U = (I\cup J)^\comp \) this recovers the usual length function of \( \tau \in S)N \). In this case we will drop the subscript \( U \). See Example \ref{ex:q-minors} for examples computing $e(\tau)$ and $\lng_U(\tau)$ for various bijections $\tau$.

\subsection{Domokos-Lenagan minors $\DLmin$ and their twists $\Tmin$, $\PTmin$}
\label{sec:DL-minors}
Let \( I \subset [N] \) be a subset of \( k \) elements, and let \( \Sym(I) \) be the symmetric group permuting the elements of \( I \).  Let us denote by \( \tau_I \) the unique order preserving bijection \( \tau_I:I \rightarrow [k] \), and for a permutation \( \sigma \in S_I \) denote \( \sigma^\circ = \tau_I \sigma \tau_I^{-1} \in S_k \).  By construction, we have
\[ \lng(\sigma^\circ) = \# \setc{(i,j) \in I\times I}{i < j \text{ and } \sigma(i) > \sigma(j)}, \]
i.e. $\lng(\sigma^\circ)$ ignores elements of $[N]$ which are not in $I$.  Recall that \( \wt(I) = \sum_{i\in I} i \). If \( I \subset [N] \) is a set with \( k \) elements, we will set \( I = \set{i_1,i_2,\ldots,i_k} \) where \( i_1 < i_2 < \ldots < i_k \).

\begin{Definition}
  \label{def:dl-minors}
For subsets \( I,J \subset [N] \) such that \( \# I = \# J = k \), we define:
\begin{enumerate}
\item The \emph{Domokos-Lenagan minor} \cite[Section 7]{DL03}
\begin{equation}
  \label{eq:DL-minor}
  \DLmin(I,J) = q^{-\wt(I)-\wt(J)}\sum_{\sigma \in \Sym(J)} (-q)^{\lng(\sigma^\circ)}x^{i_1}_{\sigma(j_1)} x^{i_2}_{\sigma(j_2)} \cdots x^{i_k}_{\sigma(j_k)}\in \FRT,
\end{equation}
where \( i_m = \tau_I^{-1}(m) \) and \( j_m = \tau_J^{-1}(m) \).
\item For \( U \subset (I \cup J)^\comp \), the \emph{truncated minor},
\begin{equation*}
  \PTmin_U(I,J) = q^{-\wt(I)-\wt(J)}\sum_{\tau:I\rightarrow J} (-q)^{\lng_U(\tau)} q^{e(\tau)} a^{i_1}_{\tau(i_1)} a^{i_2}_{\tau(i_2)} \cdots a^{i_k}_{\tau(i_k)} \in\REA,
\end{equation*}
where the sum is taken over all bijections \( I \longrightarrow J \).
\item The \emph{quantum minor},
\begin{equation*}
  \Tmin(I,J) = \Psi(\DLmin(I,J))\in\REA.
\end{equation*}
\end{enumerate}

\end{Definition}

\begin{Example}\label{ex:q-minors}
In the case \( N=4 \) and $U=(I\cup J)^\comp$, we have:
  \begin{align*}
    \DLmin(\set{1,3},\set{3,4}) &= q^{-11} \Big( x^1_3x^3_4 - q x^1_4 x^3_3 \Big), \\
    \PTmin_U(\set{1,3},\set{3,4}) &= q^{-11} \Big( -q^3 a^1_3a^3_4 + q^3 a^1_4 a^3_3 \Big),\\
    \Tmin(\set{1,3},\set{3,4}) &= q^{-11} \Big( q a^1_3 a^3_4 - q a^1_4 a^3_3 + (q-q^{-1}) a^1_4a^4_4  \Big)
  \end{align*}
  \begin{align*}
    &\DLmin(\set{1,2,4},\set{1,2,3}) \\
    &\phantom{===}= q^{-13} \Big( x^1_1 x^2_2 x^4_3 - q x^1_2 x^2_1 x^4_3 -q x^1_1 x^2_3 x^4_2 
      -q^3 x^1_3 x^2_2 x^4_1 +q^2 x^1_2 x^2_3 x^4_1 +q^2 x^1_3 x^2_1 x^4_2 \Big),\\
    &\PTmin_U(\set{1,2,4},\set{1,2,3}) \\
    &\phantom{===}= q^{-13} \Big( a^1_1 a^2_2 a^4_3 - q^2 a^1_2 a^2_1 a^4_3 -q^2 a^1_1 a^2_3 a^4_2 
      -q^4 a^1_3 a^2_2 a^4_1 +q^4 a^1_2 a^2_3 a^4_1 +q^3 a^1_3 a^2_1 a^4_2 \Big),\\
    &\Tmin(\set{1,2,4},\set{1,2,3}) \\
    &\phantom{===}= q^{-13} \Big( a^1_1 a^2_2 a^4_3 - q^2 a^1_2 a^2_1 a^4_3 -q a^1_1 a^2_3 a^4_2 
      -q^3 a^1_3 a^2_2 a^4_1 +q^3 a^1_2 a^2_3 a^4_1 +q^2 a^1_3 a^2_1 a^4_2 \\ &\phantom{==== q^{-13} \Big(}- (q-q^{-1})qa^1_3a^3_1a^4_3 + (q-q^{-1})q^2 a^1_3a^3_3a^4_1  \Big)
  \end{align*}
  The formulas above illustrate several key features.  While the Domokos-Lenagan minors and the truncated minors can be expressed as a combinatorial sum over bijections, the quantum minors cannot.  Similarly, while Propositions \ref{prp:DL-row-exp} and \ref{prp:row-exp-PTmins} below assert a row expansion formula for the Domokos-Lenagan minors and the truncated minors, the quantum minors do not admit such an obvious expansion.

In addition to extra summands appearing in the quantum minors, the coefficients of the common summands also disagree, owing to the differing notion of length appearing in their construction.  It is therefore rather remarkable that these discrepancies telescope away in the sums computed in Theorem \ref{thm:clique-sums-thm}.\end{Example}

\begin{Theorem}[\cite{DL03}, Proposition 7.2] The elements,\label{thm:DLinv}
$$\DLinv{k} = \sum_{I\in {[N] \choose k}} \DLmin(I,I),$$
are co-invariants for the adjoint $\FRT$-coaction on itself.
\end{Theorem}

\begin{Example}\label{DLinv-n=2and3}
  Let us list the Domokos-Lenagan FRT co-invariants alongside the corresponding REA invariants, when $N=2$ and $3$ for comparison:
  \begin{equation*}
    \begin{aligned}
      N&=2: \\
      \DLinv{1} &= q^{-2}x^1_1 + q^{-4} x^2_2, \qquad \DLinv{2} = q^{-6}\left(x^1_1x^2_2 - qx^1_2x^2_1\right).\\
      c_1 &= q^{-2}a^1_1 + q^{-4} a^2_2, \qquad c_2 = q^{-6}\left(a^1_1a^2_2 - q^2a^1_2a^2_1\right).\\ \\
      N&=3: \\
      \DLinv{1} &= q^{-2}x^1_1+q^{-4}x^2_2 + q^{-6}x^3_3,\\
      c_1 &= q^{-2}a^1_1+q^{-4}a^2_2 + q^{-6}a^3_3,\\
      \DLinv{2} &= q^{-6}(x^1_1x^2_2-q x^1_2x^2_1) + q^{-8}(x^1_1x^3_3 - q x^1_3x^3_1) + q^{-10} ( x^2_2x^3_3 - q x^2_3x^3_2),\\
      c_2 &= q^{-6}(a^1_1a^2_2-q^2a^1_2a^2_1) + q^{-8}(a^1_1a^3_3 - q^4a^1_3a^3_1) + q^{-10} ( a^2_2a^3_3 - q^2a^2_3a^3_2),\\
      \DLinv{3} &= q^{-12}\left(x^1_1x^2_2x^3_3 - q x^1_1x^2_3x^3_2 - q x^1_2x^2_1x^3_3 
      - q^3 x^1_3x^2_2x^3_1 + q^2 x^1_2x^2_3x^3_1 + q^2 x^1_3x^2_1x^3_2\right),\\
      c_3 &= q^{-12}\left(a^1_1a^2_2a^3_3 - q^2 a^1_1a^2_3a^3_2 -q^2 a^1_2a^2_1a^3_3  
      -q^4 a^1_3a^2_2a^3_1 + q^4 a^1_2a^2_3a^3_1 +q^3a^1_3a^2_1a^3_2\right).
    \end{aligned}
  \end{equation*}
The formulas for $\DLinv{k}$ and $c_k$ are essentially the same when $k=1$, but differ in every other case  due to the exceedance contribution and the refined notion of length for truncated minors; the simplest case where both can be seen is in comparing $\DLinv{2}$ and $c_2$ when $N=3$.  The middle minor has a different relative weighting to the others.\end{Example}

\begin{Proposition}
  \label{prp:twist-maps-dl-resh}
The twist isomorphism $\Psi$ maps each Domokos-Lenagan invariant $\DLinv{k}$ to a constant scalar multiple of Reshetikhin's canonical central elements $c_k$.
\end{Proposition}
\begin{proof}
This follows by the fact that $\DLinv{k}$ and $\Phi(c_k)$ are uniquely determined (up to a scalar) as the unique elements of degree $k$ lying in a component $\Lambda^k(V)^*\bt\Lambda^k(V)$, for the $U_q(\gl)$-bimodule action on $\FRT$.
\end{proof}

In Corollary \ref{cor:cs-are-DL-invars}, we show that the scalar multiple is one, i.e. we show that $\Psi(\DLinv{k})=c_k$.	

\begin{Remark}
While our strategy to prove Theorem \ref{main-thm} involves replacing the quantum minors with the better-behaved truncated versions, it may be interesting nevertheless to find explicit formulas for the quantum minors. For instance, the ideal formed by the quantum minors of a fixed degree $k$ will define a quantization of the variety of matrices with rank less than or equal to $k-1$.
\end{Remark}

\subsection{Row expansion formulas}
\label{sec:row-expans-form}
A useful feature of the Domokos-Lenagan minors and the truncated minors is that they admit an inductive definition via row expansion formulas, mimicking those for classical matrices.

We will employ the following notation. For a set $X$, a subset \( S \subset X \) and and element \( a \in S \) we will denote by \( S_a \) the set \( S - \set{a} \). If \( b \in X \) then \( S^b \) will denote \( S \cup \set{b} \). We will also stack this notation so that, for instance, \( S_a^{b,c} \) will denote the set \( (S - \set{a})\cup\set{b,c} \).

\begin{Proposition}
  \label{prp:DL-row-exp}
  The following row expansion formula holds,
  \begin{equation*}
    \DLmin(I,J) = \sum_{m=k}^1 (-q)^{m-1}q^{-i_1-j_m} x^{i_1}_{j_m} \DLmin(I_{i_1}, J_{j_m}).
  \end{equation*}
\end{Proposition}
\begin{proof}
  The formula follows from the fact that if \( \sigma \in \Sym(J) \) such that \( \sigma(j_1)=j_m \) then \( \lng(\sigma^\circ) = \lng(\tau^\circ)+m-1 \) where \( \tau \in \Sym(J_{j_m}) \) is the restriction of \( \sigma \). 
\end{proof}

To state the analogous row expansion formula for the truncated quantum minors, we require the following lemma, whose proof follows easily from the definitions.

\begin{Lemma}
  \label{lem:lenth-exceedance}
Let \( \tau\map{I}{J} \) be a bijection such that \( \tau(i_1) = j_m \). Let \( \tau_m \map{I_{i_1}}{J_{j_m}} \) be the bijection obtained by restricting \( \tau \). We have the following additive properties for the length and exceedance of \( \tau \),
\begin{enumerate}
\item \( e(\tau) - e(\tau_m) = \theta(j_m - i_1) \),
\item \( \lng_U(\tau) - \lng_{U}(\tau_m) = m-1 + \gamma^{IJ}_U(m), \)
\end{enumerate}
where \( \gamma_U^{IJ}(m) = \# \setc{a \in U}{i_1 < a < j_m \text{ or } j_m < a < i_1} \).
\end{Lemma}

With this in hand, we can state the following

\begin{Proposition}
  \label{prp:row-exp-PTmins}
We have the following row expansion formula for the truncated minors:
\begin{equation*}
  \PTmin_U(I, J) = \sum_{m=1}^k (-q)^{m-1+\gamma_U^{IJ}(m)} q^{\theta(j_m-i_1)-i_1-j_m} a^{i_1}_{j_m} \PTmin_{U}(I_{i_1}, J_{j_m}).
\end{equation*}
\end{Proposition}

\begin{proof}
This follows directly from Lemma~\ref{lem:lenth-exceedance}.
\end{proof}

\begin{Example}
  \label{exm:row-exp-example}
  Let \( U=\emptyset \). The row expansion formula for \( \PTmin_U(\set{1,3,4},\set{2,3,4}) \) yields:
  \begin{align*}
    \PTmin_\emptyset(\set{1,2,4},\set{1,2,3}) &= q^{-2}a^1_1\PTmin_{\emptyset}(\set{2,4},\set{2,3}) - q^{-2}a^1_2\PTmin_{\emptyset}(\set{2,4},\set{1,3}) \\
                                      &\phantom{=}+ q^{-2}a^1_3\PTmin_{\emptyset}(\set{2,4},\set{1,2})
  \end{align*}
  which is easily verified using Example~\ref{ex:q-minors} and the fact that
  \begin{align*}
    \PTmin_{\emptyset}(\set{2,4},\set{2,3}) &= q^{-11} \Big( a^2_2a^4_3 - q^2a^2_3a^4_2 \Big) \\
    \PTmin_{\emptyset}(\set{2,4},\set{1,3}) &= q^{-10} \Big( a^2_1a^4_3 - q^2a^2_3a^4_1 \Big) \\
    \PTmin_{\emptyset}(\set{2,4},\set{1,2}) &= q^{-9} \Big( a^2_1a^4_2 - q a^2_2a^4_1 \Big).
  \end{align*}
\end{Example}

\section{Expansion cliques}
\label{sec:expansion-cliques}
Modelled on the row expansion of minors in the previous section, we introduce a combinatorial device called an \emph{expansion clique}. An expansion clique $Cl_k(I,J)$ depends on two subsets $I,J\subset [N]$ of equal cardinality: roughly speaking, the expansion clique $Cl_k(I,J)$ is the set of all complementary indices $I'$ and $J'$, such that (a scalar multiple of) \( a^I_J\PTmin_U(I',J') \) appears as a summand upon performing iterated row expansion to the minor \( \PTmin_U(I\cup I', J\cup J') \). Here \( a^I_J = a^{i_1}_{j_1}\ldots a^{i_{\#I}}_{j_{\#J}} \).

\begin{Definition} Fix an integer $k \leq N$, and sets $I,J$ of size $m<k$.  The \emph{expansion clique}, $Cl_k(I,J)$, is the set,

$$Cl_k(I,J) = \left\{ (I', J') \,\,,\,\, \begin{array}{c} I' \subset \{i_m+1,\ldots N\}\\ J' \subset \{1,\ldots, N\}\end{array} \,\, \left| \,\,\begin{array}{c} I\cup I' = J\cup J',\\ |I'|=|J'|=k-m\end{array}\right. \right\}.$$
\end{Definition}

The following lemma will be useful in the proof of the main theorem.

\begin{Lemma}
  \label{lem:bijection-cliques}
  For any pair of sets \( I,J \), the map
  \begin{equation*}
    \beta : \bigsqcup_{m=1}^{k-\# I} Cl_k(I,J) \longrightarrow \bigsqcup_{s,t} Cl_k(I^s,J^t)
  \end{equation*}
  given by sending \( (I',J')_m \) to \( (I'_{i'_1},J'_{j'_m}) \in Cl_k(I^{i'_1} | J^{j'_m}) \) is a bijection. Here we label elements of the domain \( (I'J')_m \) as pairs indexed by an integer to indicate which member of the disjoint union they belong to. The \( s \) and \( t \) in the second disjoint union range over \( s \in [i_{\# I}+1,N] \) and \( t \in [N] - J \).
\end{Lemma}

\begin{proof}
The inverse is given by sending \( (I'',J'') \in Cl_k(I^s|J^t) \) to \( (I''\cup s, J''\cup t)_m \) where \( m \) is the index of \( t \) in \( J'' \cup t \).
\end{proof}

Here we collect a number of identities that will be useful later.  The reader may wish to skip ahead to Section \ref{sec:main-theorem} to see where these identities are used. Fix a pair of subsets \( I,J \subset [N] \) such that \( \# I = \# J \le k \). Fix also \( s > i_{\# I} \) and \( t \in J^\comp \), and \( (I'',J'') \in Cl_k(I^s,J^t) \) and an auxiliary set \( U \subset (I^s \cup I'')^\comp \).  We will use \( [s,t] \) to denote the interval \( [s,t] \) or \( [t,s] \) depending on whether \( s\le t \) or \( t\le s \). We use \( \ind_I(i) \) to denote the index of \( i \) in \( I \), i.e. \( \ind(i_e) = e \) and we will adopt the convention that if \( r \notin I \) then \( \ind_I(r) := \ind_{I^r}(r) \).

We shall consider the quantities:
\begin{align}
    X_{IJ}(I'',J'',s,t) &= \ind_{J''}(t) - 1 + \gamma_U^{(I'')^s, (J'')^t}(\ind_{J''}(t))     \label{eq:X-expression}\\
    Y_{IJ}(I'',J'',s,t,r) &= X_{IJ}(s,r) + \lng_{U_{sr}}(\tau_{(I'')_t^r J''}) + N((I'')_t^r,r,t),\label{eq:Y-expression}
\end{align}
where \( N(K,r,t) = \# K \cap (r,t) \).

\begin{Lemma}
  \label{lem:X-expression}
  We have the following identity:
  \[X_{IJ}(I'',J'',s,t)= \# \left( s,t \right) \cap (U \cup I \cup I'') + 1 + \ind_{I}(s) - \ind_{J}(t) - 2 \#I \cap [s,t]. \]
  In particular, it is independent of \( I'' \) and \( J'' \), so we will henceforth abbreviate $X_{IJ}(s,t)=X_{IJ}(I'',J'',s,t)$.
\end{Lemma}

\begin{proof}
  The proof will use the key fact that in the interval between \( s \) and \( t \), either only elements of \( I \) occur (and not \( I'' \)) or only elements of \( I'' \) occur (and not elements of \( I \)), depending on whether \( s \le t  \) or \( t\le s \).

  First consider the case when \( s \le t \). In this case, since every element of \( (I\cup I'')\cap[s,t] \) is actually contained in \( I'' \) (and not \( I \)), we can write
  \begin{align*}
    \ind_{J''}(t) -1 &= \# \setc{j \in J''}{j< t} \\
                         &= \# \setc{j \in J^t\cup J''}{j < t} - \#\setc{j \in J}{j<t} \\
                         &= \# \setc{i \in I^s\cup I''}{i < t} - \#\setc{j \in J}{j<t} \\
                         &= \# \setc{i \in I^s}{i \le s} + \#\setc{i \in I''}{s<i<t} - \#\setc{j \in J}{j<t} \\
                         &= \ind_{I}(s) - \ind_{J}(t)+1 + \#\setc{i \in I''}{s<i<t}.
  \end{align*}
  On the other hand, note that the smallest element of \( (I'')^s \) is \( s \) so
  \begin{align*}
    \gamma_U^{(I'')^s,(J'')^t}(\ind_{J''}(t)) &= \# (s,t)\cap U.
  \end{align*}
  Taking the sum gives the desired result (noting that \( \#I \cap [s,t] = 0 \) when \( s \le t \)).

  In the case \( t \le s \), we proceed similarly. So
  \begin{align*}
    \ind_{J''}(t) -1 &= \# \setc{i \in I^s\cup I''}{i < t} - \#\setc{j \in J}{j<t} \\
                        &= \# \setc{i \in I^s}{i \le s} - \#\setc{i \in I}{t \le i \le s} - \#\setc{j \in J}{j<t} \\
                         &= \ind_{I}(s) - \ind_{J}(t)+1 - \#\setc{i \in I}{t \le i \le s},
  \end{align*}
  and again \( \gamma_U^{(I'')^s,(J'')^t}(\ind_{J''}(t)) = \# (s,t)\cap U  \) so taking the sum again gives the desired result.
\end{proof}

For any two subsets \( I,J \subset [N] \) such that \( \# I = \# J \), let \( \tau_{IJ} \) be the unique, order preserving bijection \( I \longrightarrow J \).

\begin{Lemma}
  \label{lem:telescoping-condition}
  Assume that \( s < t \). If \( r < r' \) are adjacent elements of \( (s,t] - I'' \) then
  \[ Y_{IJ}(I'',J'',s,t,r') - Y_{IJ}(I'',J'',s,t,r) = 2. \]
\end{Lemma}

\begin{proof}
  We can use the expression from Lemma~\ref{lem:X-expression} to see that
  \begin{align*}
    X_{IJ}(s,r')-X_{IJ}(s,r) &= \#(s,r')\cap (U \cup I \cup I'') +\ind_{I}(s)-\ind_{J}(r') \\
                             &\phantom{=}- \left( \#(s,r)\cap (U \cup I \cup I'')+\ind_{I}(s)-\ind_{J}(r) \right) \\
                             &= \#[r,r')\cap (U \cup I \cup I'') - \left( \ind_{J}(r') - \ind_{J}(r) \right) \\
                             &= \# I'' \cap [r,r'] +1 - \# J \cap [r,r'].
  \end{align*}
  The last equality is due to the fact that \( r \) and \( r' \) are adjacent in \( (s,t] - I'' \), thus every element of \( (r,r') \) must be in \( I'' \). Clearly
  \[ N(I'',r',t) - N(I'',r,t) = -N(I'',r,r') = -\#I'' \cap [r,r']. \]

  Finally, we will show that
  \[ \lng_{U}(\tau_{(I'')_t^{r'} J''}) - \lng_{U}(\tau_{(I'')_t^r J''}) = \# J \cap [r,r'] + 1 \]
  which will mean we are done. To do this, we imagine the sets \( I,I'' \) depicted by circles and crosses (similarly the sets \( J, J'' \)) and strings connecting \( I'' \) and \( J'' \) to depict \( \tau_{(I'')^r_t, J''} \) as in Figure~\ref{fig:bijection--tauUrt}. We also depict with red lines, the elements of \( U \). We call these \emph{loose strings}. 
  
  \begin{figure}
    \centering
    \begin{tikzpicture}[inner sep=0.5pt]
      \node (i1) at (1,2) {\( \circ \)};
      \node (i2) at (2,2) {\( \circ \)};
      \node (i3) at (3,2) {\( \circ \)};
      \node (i4) at (4,2) {\( \circ \)};
      \node (i5) at (5,2) {\( \circ \)};
      \node (s) at (6,2) {\( \star \)};
      \node (ip1) at (7,2) {\( \times \)};
      \node (ip2) at (8,2) {\( \times \)};
      \node (ip3) at (9,2) {\( \times \)};
      \node (ip4) at (10,2) {\( \times \)};
      \node (ip5) at (11,2) {\( \times \)};
      \node (ip6) at (12,2) {\( \times \)};

      \node (jp1) at (1,0) {\( \times \)};
      \node (j1) at (2,0) {\( \circ \)};
      \node (jp2) at (3,0) {\( \times \)};
      \node (jp3) at (4,0) {\( \times \)};
      \node (j2) at (5,0) {\( \circ \)};
      \node (j3) at (6,0) {\( \circ \)};
      \node (j4) at (7,0) {\( \circ \)};
      \node (jp4) at (8,0) {\( \times \)};
      \node (jp5) at (9,0) {\( \times \)};
      \node (t) at (10,0) {\( \star \)};
      \node (jp6) at (11,0) {\( \times \)};
      \node (j5) at (12,0) {\( \circ \)};

      \node[anchor=south] at (6,2.2) {\( s \)};
      \node[anchor=south] at (10,2.2) {\( t \)};
      \node[anchor=south] at (7.67,2.2) {\( r \)};
      \node[anchor=south] at (9.5,2.2) {\( r' \)};

      \node (u1) at (1.5,2) {\( \cdot \)};
      \node (l1) at (1.5,0) {\( \cdot \)};
      \node (u2) at (3.5,2) {\( \cdot \)};
      \node (l2) at (3.5,0) {\( \cdot \)};
      \node (u3) at (5.33,2) {\( \cdot \)};
      \node (l3) at (5.33,0) {\( \cdot \)};
      \node (u4) at (5.67,2) {\( \cdot \)};
      \node (l4) at (5.67,0) {\( \cdot \)};
      \node (u5) at (6.33,2) {\( \cdot \)};
      \node (l5) at (6.33,0) {\( \cdot \)};
      \node (u6) at (7.33,2) {\( \cdot \)};
      \node (l6) at (7.33,0) {\( \cdot \)};
      \node (u7) at (7.67,2) {\( \cdot \)};
      \node (l7) at (7.67,0) {\( \cdot \)};
      \node (u8) at (9.5,2) {\( \cdot \)};
      \node (l8) at (9.5,0) {\( \cdot \)};
      \node (u9) at (10.5,2) {\( \cdot \)};
      \node (l9) at (10.5,0) {\( \cdot \)};

      \draw (ip1) to (jp1);
      \draw (u7) to (jp2);
      \draw (ip2) to (jp3);
      \draw (ip3) to (jp4);
      \draw (ip5) to (jp5);
      \draw (ip6) to (jp6);

      \draw[red] (u1) to (l1);
      \draw[red] (u2) to (l2);
      \draw[red] (u3) to (l3);
      \draw[red] (u4) to (l4);
      \draw[red] (u5) to (l5);
      \draw[red] (u6) to (l6);
      \draw[red,dashed] (u7) to (l7);
      \draw[red] (u8) to (l8);
      \draw[red] (u9) to (l9);

      \draw[red] (ip4) to (t);
    \end{tikzpicture}
    \caption{The bijection \( \tau_{(I'')^r_t, J''} \). \( \circ \in I \) or \( J \), and \( \times \in I'' \) or \( J'' \).}
    \label{fig:bijection--tauUrt}
  \end{figure}
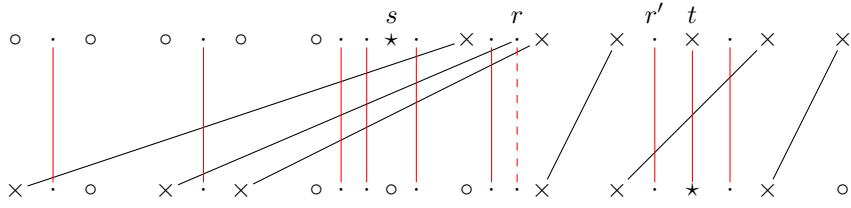

  Since \( \tau_{(I'')^r_t,J''} \) is order preserving, its length is exactly the number of times black strings cross loose strings. In the example in Figure~\ref{fig:bijection--tauUrt}, the length is \( 16 \). We can imagine changing the diagram to the diagram for the bijection \( \tau_{(I'')^{r'}_t,J''} \) by moving the top of the \( r \)-string and the top of each \( I''\cap[r,r'] \)-sting, up one available position (with the last string moving to \( r' \)). We only need to account for crossings gained and lost.

  Since \( r \) and \( r' \) are adjacent in \( (s,t] - I'' \), these considerations simplify considerably: crossings are only lost by removing \( r' \) as a loose string and crossings are only gained by adding \( r \) as a loose string. The number of strings crossing \( r \) (i.e. the number of crossings gained) is
  \[ \ind_{J''}(r)  - \ind_{I''}(r) \]
  and the number of strings crossing \( r' \) (crossings lost) is
  \[ \ind_{J''}(r') - 1 - \ind_{I''}(r') . \]
  The difference of these two quantities is
  \[ \# I'' \cap (r,r'] - \# J'' \cap (r,r'] + 1, \]
  so all that is left to notice is that \( \# I'' \cap (r,r'] - \# J'' \cap (r,r'] = \# J \cap [r,r']  \).
\end{proof}

\begin{Lemma}
  \label{lem:Y-min-r}
  Assume that \( s < t \). Let \( r \) be the minimal element of \( (s,t]-I'' \). Then
  \begin{equation*}
    Y_{IJ}(I'',J'',s,t,r) = \lng_{U}(\tau_{(I'')^s (J'')^t}) + \ind_{J''}(t)-1.
  \end{equation*}
\end{Lemma}

\begin{proof}
Recall the definitions \eqref{eq:X-expression}, \eqref{eq:Y-expression}:
  \begin{align*}
    X_{IJ}(s,r) &= \ind_{J''}(r) - 1 + \gamma_U^{(I'')^s,(J'')^r}(\ind_{J''}(r)),\\
    Y_{IJ}(I'',J'',s,t,r) &= X_{IJ}(s,r) + \lng_{U_{sr}}(\tau_{(I'')^r_t J''}) + N(I'',r,t).
  \end{align*}
  Note that since \( r \) is chosen to be minimal, \( \gamma_U^{(I'')^s,(J'')^r}(\ind_{J''}(r)) = 0 \). Thus the claim in the Lemma reduces to proving
  \begin{align*}
    \lng_{U}(\tau_{(I'')^s (J'')^t}) - \lng_{U}(\tau_{(I'')^r_t J''}) &= N(I'',r,t) + \ind_{J''}(r) - \ind_{J''}(t)\\
    & = N(I'',r,t) - N(J'',r,t).
  \end{align*}

  We will start with the diagram for \( \tau_{(I'')^r_t J''} \) (see the example in Figure~\ref{fig:bijection--tau-r-min}) and describe a process by which to transform it into the diagram for \( \tau_{(I'')^s (J'')^t} \). At each step we will account for any crossings gained or lost and thus arrive at an expression for the difference of the lengths.

  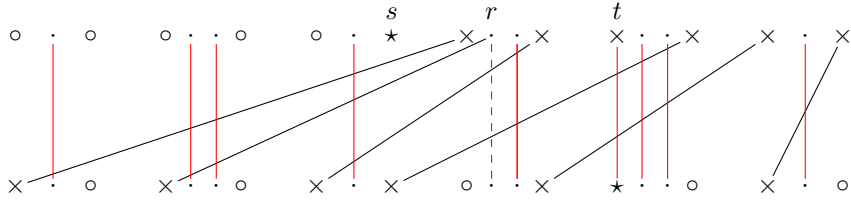
\begin{figure}
    \centering
    \begin{tikzpicture}[inner sep=0.5pt]
      \node (i1) at (1,2) {\( \circ \)};
      \node (i2) at (2,2) {\( \circ \)};
      \node (i3) at (3,2) {\( \circ \)};
      \node (i4) at (4,2) {\( \circ \)};
      \node (i5) at (5,2) {\( \circ \)};
      \node (s) at (6,2) {\( \star \)};
      \node (ip1) at (7,2) {\( \times \)};
      \node (ip2) at (8,2) {\( \times \)};
      \node (ip3) at (9,2) {\( \times \)};
      \node (ip4) at (10,2) {\( \times \)};
      \node (ip5) at (11,2) {\( \times \)};
      \node (ip6) at (12,2) {\( \times \)};

      \node (jp1) at (1,0) {\( \times \)};
      \node (j1) at (2,0) {\( \circ \)};
      \node (jp2) at (3,0) {\( \times \)};
      \node (j2) at (4,0) {\( \circ \)};
      \node (jp3) at (5,0) {\( \times \)};
      \node (jp4) at (6,0) {\( \times \)};
      \node (j3) at (7,0) {\( \circ \)};
      \node (jp5) at (8,0) {\( \times \)};
      \node (t) at (9,0) {\( \star \)};
      \node (j4) at (10,0) {\( \circ \)};
      \node (jp6) at (11,0) {\( \times \)};
      \node (j5) at (12,0) {\( \circ \)};

      \node[anchor=south] at (6,2.2) {\( s \)};
      \node[anchor=south] at (9,2.2) {\( t \)};
      \node[anchor=south] at (7.33,2.2) {\( r \)};

      \node (u1) at (1.5,2) {\( \cdot \)};
      \node (l1) at (1.5,0) {\( \cdot \)};
      \node (u2) at (3.33,2) {\( \cdot \)};
      \node (l2) at (3.33,0) {\( \cdot \)};
      \node (u3) at (3.67,2) {\( \cdot \)};
      \node (l3) at (3.67,0) {\( \cdot \)};
      \node (u4) at (5.5,2) {\( \cdot \)};
      \node (l4) at (5.5,0) {\( \cdot \)};
      \node (u5) at (7.33,2) {\( \cdot \)};
      \node (l5) at (7.33,0) {\( \cdot \)};
      \node (u6) at (7.67,2) {\( \cdot \)};
      \node (l6) at (7.67,0) {\( \cdot \)};
      \node (u8) at (9.33,2) {\( \cdot \)};
      \node (l8) at (9.33,0) {\( \cdot \)};
      \node (u9) at (9.67,2) {\( \cdot \)};
      \node (l9) at (9.67,0) {\( \cdot \)};
      \node (u10) at (11.5,2) {\( \cdot \)};
      \node (l10) at (11.5,0) {\( \cdot \)};

      \draw (ip1) to (jp1);
      \draw (u5) to (jp2);
      \draw (ip2) to (jp3);
      \draw (ip4) to (jp4);
      \draw (ip5) to (jp5);
      \draw (ip6) to (jp6);

      \draw[red] (u1) to (l1);
      \draw[red] (u2) to (l2);
      \draw[red] (u3) to (l3);
      \draw[red] (u4) to (l4);
      \draw[red,dashed] (u5) to (l5);
      \draw[red] (u6) to (l6);
      \draw[red] (u7) to (l7);
      \draw[red] (u8) to (l8);
      \draw[red] (u9) to (l9);
      \draw[red] (u10) to (l10);

      \draw[red] (ip3) to (t);
    \end{tikzpicture}
    \caption{The bijection \( \tau_{(I'')^r_t, J''} \) for \( r \) minimal.}
    \label{fig:bijection--tau-r-min}
  \end{figure}

  The first step is to move the strings starting at positions \( i''_1, i''_2, \ldots, r \in (I'')^r \) one position to the left. Since \( r \) is minimal, no crossings are lost, and no new crossings are gained.

  The second step will be to add in the loose string at position \( r \). This creates as many new crossings as there are strings crossing \( r \), i.e. there will be
  \begin{equation}
    \label{eq:gains-1}
    \ind_{J''}(r) - 1 - \ind_{I''}(s) 
  \end{equation}
  new crossings.

  Thirdly, if \( m = \ind_{J''}(t) \) and \( e = \ind_{I''}(t) \), we take the \( i''_e \)-string, the \( i''_{e+1} \)-string, up to the \( i''_{m-1} \)-string, and move them all one position to the left. This reduces the number of crossings by
  \begin{equation}
    \label{eq:losses-1}
    \sum_{k=e}^{m-1} \# U \cap (i''_{k-1},i''_k) = \# U \cap (t,i''_{m-1}) - (m-e-1). 
  \end{equation}

  The fourth step will be to introduce a string starting at \( i''_{m-1} \) and ending at \( t \). This introduces crossings equal to the number of loose strings between \( t \) and \( i''_{m-1} \), i.e. exactly the number given in~(\ref{eq:losses-1}). Thus these contributions cancel out.

  Lastly, we delete the loose string at position \( t \). This results in a loss of crossings equal to the number of strings crossing \( t \), i.e. there will be
  \begin{equation}
    \label{eq:losses-2}
    \ind_{J''}(t)-1 - \ind_{I''}(t)
  \end{equation}
  crossings gained.

  The total number of crossings gained is thus given by the difference of (\ref{eq:gains-1}) and (\ref{eq:losses-2}). Noting that
  \begin{align*}
    N(I'',r,t) &= \ind_{I''}(t) - \ind_{I''}(r), \text{ and } \\
    N(J'',r,t) &= \ind_{J''}(t) - \ind_{J''}(r), 
  \end{align*}
  the claim is proved.
\end{proof}

\section{The formula for the invariants $c_k$}
\label{sec:main-theorem}
In this section we prove our main results, by an elaborate induction using the row expansion formulas for the Domokos-Lenagan and truncated minors.  They key observation is that while the truncated minors only approximate the image under the twisting isomorphism $\Psi$, the symmetric sums appearing in the definition of the invariants coincide whether we use the truncated or the quantum minors.  This is captured in the following theorem.

\begin{Theorem}
  \label{thm:clique-sums-thm}
  For any \( 1 \le k \le N \), and any pair \( I,J \subset \brak{N} \). We have
  \begin{equation}
    \label{eq:main-identity}
    \sum_{I',J'} \PTmin_{(I\cup I')^\comp}(I',J') = \sum_{I',J'} (-q)^{\lng_{(I\cup I')^\comp}(\tau_{I'J'})} \Tmin(I',J'),
  \end{equation}
where both sums are over \( (I',J') \in Cl_k(I,J) \).
\end{Theorem}

\begin{proof}
  Recall that \( \Tmin(I',J') = \Psi\left( \DLmin(I',J') \right) \). We will prove the identity by inducting on \( k-m \). When \( k=m \), the clique \( Cl_k(I,J) \) is either empty and the result is vacuous, or \( I=J \) and the clique contains only the pair of empty subsets \( (\emptyset, \emptyset) \). In this case \( \DLmin(\emptyset,\emptyset) = 1 \) and so \( \PTmin_{I^\comp}(\emptyset,\emptyset) = 1 = \Tmin(\emptyset,\emptyset) \). When \( k > m \), we will in fact prove
  \[ \sum_{I',J'} \Phi(\PTmin_{(I\cup I')^\comp}(I',J')) = \sum_{I',J'} (-q)^{\lng_{(I\cup I')^\comp}(\tau_{I'J'})} \DLmin(I',J'). \]
  Since \( \Phi \) is an isomorphism (of vector spaces), this is equivalent to~(\ref{eq:main-identity}).  However, the PBW ordering relations in the FRT algebra are much simpler, so it is convenient to check the identity there.

To begin, we use the row expansion formula from Proposition~\ref{prp:row-exp-PTmins} to express the left hand side as
  \begin{align*}
  \sum_{I',J'} \Phi&\left(\PTmin_{(I\cup I')^\comp}(I',J')\right) \\ &\phantom{=}= \sum_{I', J'} \sum_{m=1}^{k-\# I} q^{\theta(j_m'-i_1')-i'_1-j'_m}(-q)^{m-1+\gamma_{(I\cup I')^\comp}^{I'J'}(m)} \Phi\left(a^{i'_1}_{j'_m} \PTmin_{(I\cup I')^\comp}(I'_{i'_1}, J'_{j'_m})\right).
  \end{align*}
  Both sums are taken over \( I',J' \in Cl_k(I,J) \).  The terms in the inner summation on the RHS are indexed by the domain of the bijection \( \beta \) in Lemma~\ref{lem:bijection-cliques}. We can use this bijection to re-index the sum and obtain the following expression for the RHS above:
  \begin{align*}
&\sum_{I',J'} \Phi\left(\PTmin_{(I\cup I')^\comp}(I',J')\right) \\ &\phantom{}= \sum_{s,t}\sum_{I'', J''} q^{- s - t + \theta(t-s)} (-q)^{\ind_{J''}(t) -1 + \gamma_{(I^s\cup I'')^\comp}^{(I'')^s, (J'')^t}(\ind_{J''}(t))} \Phi\left( a^{s}_{t} \PTmin_{{(I^s\cup I'')^\comp}}(I'', J'') \right).
   \end{align*}
   The first sum is taken over all pairs \( s,t \) so that \( s \in [i_{\# I}+1, N] \) and \( t \in [N] - J \). Here we recognize the exponent of $(-q)$ as the quantity,
\[ X_{IJ}(s,t) = \ind_{J''}(t) - 1 + \gamma_{(I^s\cup I'')^\comp}^{(I'')^s,(J'')^t}(\ind_{J''}(t)), \]
which was defined in \eqref{eq:X-expression}.  {\it A priori} this quantity depends on \( I'',J'' \), however Lemma~\ref{lem:X-expression} tells us this is independent of \( I'' \) and \( J'' \). Thus we can take this factor and the \( a^s_t \) out of the first summand to obtain
   \begin{align*}
     \sum_{I',J'} \Phi&\left(\PTmin_{(I\cup I')^\comp}(I',J')\right) \\ &\phantom{=======}=  \sum_{s,t}q^{-s-t+\theta(t-s)}(-q)^{X_{IJ}(s,t)} \Phi\left(  a^s_t \sum_{I'', J''}  \PTmin_{(I^s\cup I'')^\comp}(I'', J'') \right).
   \end{align*}
   We can now apply the induction hypothesis to obtain
   \begin{align*}
     \sum_{I',J'} \Phi&\left(\PTmin_{(I\cup I')^\comp}(I',J')\right) \\
     &\phantom{=}= \sum_{s,t}q^{-s-t+\theta(t-s)}(-q)^{X_{IJ}(s,t)} \Phi\left( a^s_t \sum_{I'', J''} (-q)^{\lng_{(I^s\cup I'')^\comp}(\tau_{I'',J''})} \Tmin(I'', J'') \right) \\
     &\phantom{=}= \sum_{s,t}q^{-s-t+\theta(t-s)}(-q)^{X_{IJ}(s,t)} \sum_{I'', J''} (-q)^{\lng_{(I\cup I')^\comp}(\tau_{I'',J''})} \Phi\left( a^s_t\Tmin(I'', J'') \right).
   \end{align*}
We now require the following:

\begin{Lemma}\label{lem-aijtimesDL} Let \( I,J \subset [N] \) be subsets of the same size. Assume that \( i \in [N] \) is smaller than any element of \( I \). If $j\not\in I$, we have $\Phi(a^i_j\Tmin(I,J)) = x^i_j\DLmin(I,J)$.  If $j\in I$ we have:
  \begin{align*}\Phi(a^i_j\Tmin(I,J)) &= q^{-1}x^i_j\DLmin(I,J) + (q^{-1}-q)\sum_{k>j} q^{j-k}(-q)^{N(I,j,k)}x^i_k\DLmin(I^k_j,J)\end{align*}
  where \( N(I,j,k) = \#I \cap (j,k) \).
\end{Lemma}
\begin{proof}
  If \( j \notin I \), then only the first case in Lemma~\ref{lem:expanded-rea-mult} applies and the formula \( \Phi(a^i_j\Tmin(I,J)) = x^i_j \DLmin(I,J) \) follows. If \( j \in I \) then a combination of cases (1) and (3) in Lemma~\ref{lem:expanded-rea-mult} imply the formula. 
\end{proof}

Applying the calculation from Lemma~\ref{lem-aijtimesDL} we arrive at the expression
   \begin{align*}
     \sum_{I',J'} \Phi&\left(\PTmin_{(I\cup I')^\comp}(I',J')\right)  \\
     &\phantom{==}= \sum_{s \ge t} q^{-s-t+\theta(t-s)}(-q)^{X_{IJ}(s,t)} \sum_{I'', J''} (-q)^{\lng_{{(I^s\cup I'')^\comp}}(\tau_{I'',J''})}  x^{s}_{t} \DLmin(I'', J'')  \\
  &\phantom{===}+ \sum_{s < t} q^{-s-t+\theta(t-s)}(-q)^{X_{IJ}(s,t)} \sum_{I'', J''} (-q)^{\lng_{{(I^s\cup I'')^\comp}}(\tau_{I'',J''})} \\
  &\phantom{===}\cdot \left( q^{-1} x^{s}_{t}  \DLmin(I'', J'') + (q^{-1}-q) \sum_{r>t} q^{t-r}(-q)^{N(I'',t,r)} x^s_r \DLmin((I'')^r_t,J'')  \right). 
   \end{align*}
   The next step will be to swap the roles of \( r \) and \( t \) in order to reorganize the sum. The aim is to collect together the terms \( x^s_t \DLmin(I'',J'') \). The resulting expression is
   \begin{align*}
     \sum_{I',J'} \Phi\left(\PTmin_{(I\cup I')^\comp}(I',J')\right) &=\sum_{s \ge t} \sum_{I'', J''} q^{-s - t} (-q)^{X_{IJ}(s,t) + \lng_{(I^s\cup I'')^\comp}(\tau_{I'',J''})}  x^{s}_{t} \DLmin(I'', J'')  \\
  + &\sum_{s < t} \sum_{I'', J''} \bigg[ q^{-s - t} (-q)^{X_{IJ}(s,t)+\lng_{(I^s\cup I'')^\comp}(\tau_{I'',J''})}x^{s}_{t}  \DLmin(I'', J'') \\
  &+ (1-q^2) \sum_{\substack{r \in (s,t) \\ r \notin I''}} q^{-s - t)} (-q)^{Y_{IJ}(I'',J'',s,t,r)}x^{s}_{t}  \DLmin(I'', J'') \bigg], 
   \end{align*}
   where we recall from \eqref{eq:Y-expression} that
   \begin{equation*}
     Y_{IJ}(I'',J'',s,t,r) = X_{IJ}(s,r)+\lng_{(I^s\cup I'')^\comp}(\tau_{(I'')^r_t,J''}) + N((I'')^r_t,r,t).
   \end{equation*}

   Lemma~\ref{lem:telescoping-condition} shows that the sum indexed by \( r \) telescopes and we arrive at the expression
   \begin{align*}
\sum_{I',J'} \Phi\left(\PTmin_{(I\cup I')^\comp}(I',J')\right) &=  \sum_{s \ge t} \sum_{I'', J''} q^{-s - t} (-q)^{X_{IJ}(s,t)+ \lng_{(I^s\cup I'')^\comp}(\tau_{I'',J''})}  x^{s}_{t} \DLmin(I'', J'')  \\
  + &\sum_{s < t} \sum_{I'', J''} q^{-s - t} (-q)^{Y_{IJ}(I'',J'',s,t,r)}x^{s}_{t}  \DLmin(I'', J'') 
   \end{align*}
   where \( r \) is now the minimal element of \( (s,t] - I'' \). By the definition of \( \gamma \), as long as \( s \ge t \), we have
   \[ \lng_{(I^s\cup I'')^\comp}(\tau_{(I'')^s(J'')^t}) = \lng_{{(I^s\cup I'')^\comp}}(\tau_{I'' J''}) + \gamma_{(I^s\cup I'')^\comp}^{(I'')^s (J'')^t}(\ind_{J''}(t)) \]
   Thus using Lemma~\ref{lem:Y-min-r}, the above expression simplifies to
   \begin{align*}
\sum_{I',J'} \Phi&\left(\PTmin_{(I\cup I')^\comp}(I',J')\right) \\ &\phantom{=======}= \sum_{s,t} \sum_{I'', J''}  q^{-s -t} (-q)^{\lng_{(I^s\cup I'')^\comp}(\tau_{(I'')^s,(J'')^t})+\ind_{J''}(t)-1}  x^{s}_{t} \DLmin(I'', J'').
   \end{align*}
   This allows use to use the bijection in Lemma~\ref{lem:bijection-cliques} in reverse and obtain
   \begin{equation*}
     \sum_{I',J'} \Phi\left(\PTmin_{(I\cup I')^\comp}(I',J')\right) =\sum_{m=1}^{k-\# I}\sum_{I', J'} q^{-i'_1-j'_m}(-q)^{\lng_{(I\cup I')^\comp}(\tau_{I' J'})+m-1} x^{i'_1}_{j'_m} \DLmin(I'_{i'_1},J'_{j'_m})
   \end{equation*}
   which we recognise as the row expansion formula from Proposition~\ref{prp:DL-row-exp} applied to
   \begin{equation*}
     \sum_{I',J'} (-q)^{\lng_{(I\cup I')^\comp}(\tau_{I'J'})} \DLmin(I',J'), 
   \end{equation*}
and the proof is complete.
 \end{proof}

Our claimed formula for Reshetikhin's invariants $c_k$ is now a special case of the above, when $I=J=\emptyset$:

 \begin{Corollary}
   \label{cor:cs-are-DL-invars}
   Reshetikhin's elements \( c_k \) are given by the PBW ordered expressions
   \begin{equation*}
     c_k = \sum_{I \in {\SetN \choose k}} \PTmin_{I^\comp}(I,I) = \sum_{I \in {\SetN \choose k}}q^{-2\wt(I)}\!\!\sum_{\sigma\in Sym(I)} (-q)^{\ell(\sigma)}\cdot q^{e(\sigma)}\cdot a^{i_1}_{\sigma(i_1)}\cdots a^{i_k}_{\sigma(i_k)}.
   \end{equation*}
 \end{Corollary}

 \begin{proof}
   By Proposition~\ref{prp:twist-maps-dl-resh}, $c_k$ is a scalar multiple of \( \Psi(\DLinv{k}) = \sum_{I} \Tmin(I,I) \). We apply Theorem~\ref{thm:clique-sums-thm} with \( I=J=\emptyset \). Then \( Cl_k(I,J) \) is the set of all pairs \( (I',I') \) where \( I' \) is a \( k \)-element subset of \( [N] \). Thus the left hand side of \ref{eq:main-identity} is \( \sum_{I} \PTmin_{I^\comp}(I,I) \). Note that \( \tau_{I'I'}=\id \) so the right hand side of (\ref{eq:main-identity}) is precisely \( \Psi(\DLinv{k}) \). All that remains is to check the scalar, and for this we can compare a single monomial.  It is clear from Lemma \ref{lem:expanded-rea-mult} that $\Psi(x^1_1\cdots x^k_k)=a^1_1\cdots a^k_k$, hence the required scalar multiple must be one.
 \end{proof}

\printbibliography

\end{document}